\documentclass[10pt,francais,a4paper]{smfart}

\usepackage{smfthm,smfenum}
\usepackage[frenchb]{babel}
\usepackage[latin1]{inputenc}
\usepackage[T1]{fontenc}
\usepackage{t1enc}
\usepackage{xy}
\usepackage{eepic,epic}
\usepackage{epsf,psfig} 
\usepackage{mathrsfs}
\usepackage{graphics,psfrag}
\usepackage{graphicx} 
\xyoption{all}

\def\Ddots{\mathinner{\mkern1mu\raise1pt\hbox{.}\mkern2mu\raise4pt\hbox{.}\mkern
2mu
  \raise7pt\vbox{\kern7pt{\hbox{.}}}\mkern1mu}}

\title{K-th{\'e}orie {\'e}quivariante des
  tours de Bott. Application {\`a} la structure multiplicative de la
   K-th{\'e}orie {\'e}quivariante des vari{\'e}tes de drapeaux}

 \alttitle{Equivariant K-theory of Bott
  towers. Application to the multiplicative structure of the equivariant 
  K-theory of flag varieties}

\author[Matthieu Willems ]{Matthieu Willems}

\address{University of Toronto \\
Department of Mathematics \\ 
100 St George Street \\
Toronto Ontario \\
Canada M5S 3G3  }

\email{matthieu.willems@polytechnique.org}

\begin{document}
    
\bibliographystyle{smfplain}

\subjclass{$19$L$47$, $14$M$15$ $14$M$25$}


 \altkeywords{equivariant K-theory, flag varieties, toric varieties}
\mainmatter
\maketitle 

\tableofcontents

\begin{abstract}

On construit une base de la $K$-th{\'e}orie {\'e}quivariante des tours de
Bott, et on d{\'e}crit pr{\'e}cis{\'e}ment la structure multiplicative de ces
alg{\`e}bres. On en d{\'e}duit des r{\'e}sultats analogues pour les
vari{\'e}t{\'e}s de Bott-Samelson. Le lien entre les vari{\'e}t{\'e}s de
drapeaux et les vari{\'e}t{\'e}s de Bott-Samelson nous permet alors de
donner une m{\'e}thode de calcul des constantes de structure de la
$K$-th{\'e}orie {\'e}quivariante des vari{\'e}t{\'e}s de drapeaux par rapport
{\`a} la base construite par Kostant et Kumar dans~\cite{kkk}.

\end{abstract}

\begin{altabstract}

We construct a basis of the equivariant $K$-theory of Bott towers, and
we describe precisely the multiplicative structure of these
algebras. We deduce similar results for Bott-Samelson
varieties. Thanks to the link between flag varieties and Bott-Samelson
varieties, we give a method to compute the structure constants of the
equivariant $K$-theory of flag varieties  in the basis constructed by
Kostant and Kumar in~\cite{kkk}.

\end{altabstract}

\section{Introduction}

Soit $G$ un groupe de Lie semi-simple complexe connexe. Soient $B
\subset G$ un sous-groupe de Borel de $G$, et $T \subset B$ un
tore compact maximal de $B$. On note $R[T]$ l'anneau des
repr{\'e}sentations de $T$ et $X=G/B$ la vari{\'e}t{\'e} de drapeaux
associ{\'e}e {\`a} ces donn{\'e}es (plus g{\'e}n{\'e}ralement, on s'int{\'e}ressera aux
vari{\'e}t{\'e}s de drapeaux  des groupes de Kac-Moody). La multiplication {\`a}
gauche dans $G$ induit une action de $T$ sur $X$. La $K$-th{\'e}orie
$T$-{\'e}quivariante de $X$, not{\'e}e $K_T(X)$, a {\'e}t{\'e} initialement
identifi{\'e}e par Kostant et Kumar dans~\cite{kkk}. Ils
construisent notamment une base $\{\hat{\psi}^w\}_{w\in W}$
de $K_T(X)$ en tant que $R[T]$-module, o{\`u} $W$ est le groupe de 
Weyl de $X$. Une fois qu'on a
construit une base de $K_T(X)$, un des probl{\`e}mes
fondamentaux est de trouver des formules pour multiplier deux
{\'e}l{\'e}ments de cette base. La m{\^e}me question se pose pour la
$K$-th{\'e}orie ordinaire de $X$, not{\'e}e $K(X)$. Dans \cite{pittie}, Pittie et
Ram  donnent une formule pour multiplier dans $K(X)$ la classe
d'un fibr{\'e} en droites par une classe $[\mathcal{O}_{\overline{X}_w}]$, o{\`u} 
pour $w \in W$, $[\mathcal{O}_{\overline{X}_w}] \in K(X)$
d{\'e}signe la classe du faisceau structural de la vari{\'e}t{\'e} de
Schubert $\overline{X}_w \subset X$. Ce r{\'e}sultat a {\'e}t{\'e}
initialement formul{\'e} pour $G=SL(n,\mathbb{C})$ par Fulton et Lascoux
dans \cite{ful-las}. Dans~\cite{+brion}, Brion d{\'e}termine
le signe des constantes de structure de $K(X)$ par rapport {\`a} la base
$\{[\mathcal{O}_{\overline{X}_w} ]\}_{w \in W}$. Un des objectifs de
cet article est de donner une m{\'e}thode de calcul g{\'e}n{\'e}rale des
constantes de structure de $K_T(X)$ par rapport {\`a} la base
$\{\hat{\psi}^w\}_{w \in W}$.

\medskip

Les vari{\'e}t{\'e}s de Bott-Samelson ont {\'e}t{\'e} initialement introduites
par Bott et Samelson dans~\cite{bs}. Ces vari{\'e}t{\'e}s sont munies
d'une action de $T$ et d'une application $T$-{\'e}quivariante {\`a}
valeurs dans $X$. Dans~\cite{demazure}, Demazure montre qu'elles 
permettent de d{\'e}singulariser les vari{\'e}t{\'e}s de Schubert.

Les tours de Bott sont des vari{\'e}t{\'e}s toriques particuli{\`e}res construites
par fibrations successives de fibre $\mathbb{C}P^1$, et sont munies
de l'action d'un tore compact $D$ de m{\^e}me dimension que la
vari{\'e}t{\'e}. Une vari{\'e}t{\'e} de Bott-Samelson $\Gamma$ munie de
l'action  du tore $T$ peut {\^e}tre vue comme une tour de Bott $Y$, et  
l'action de $T$ sur $\Gamma$ s'identifie {\`a} celle d'un sous-tore de $D$ sur $Y$.

\medskip

Dans \cite{mw3}, on a donn{\'e} une premi{\`e}re description de la
cohomologie et de la $K$-th{\'e}orie {\'e}quivariantes des vari{\'e}t{\'e}s de
Bott-Samelson, et on a explicit{\'e} le lien avec  les vari{\'e}t{\'e}s de
drapeaux. Gr{\^a}ce {\`a} ces r{\'e}sultats, on a donn{\'e} une preuve
g{\'e}om{\'e}trique des formules de restrictions aux points fixes d'une
base de la cohomologie $T$-{\'e}quivariante de $X$. Ces formules ont
{\'e}t{\'e} initialement d{\'e}montr{\'e}es par Sara Billey dans~\cite{sb}. De
plus, on a trouv{\'e} une formule pour les restrictions aux points fixes
de la base $\{\hat{\psi}^w\}_{w \in W}$ de $K_T(X)$. Ces formules ont {\'e}t{\'e}
d{\'e}t{\'e}rmin{\'e}es ind{\'e}pendamment par William Graham gr{\^a}ce {\`a} d'autres
m{\'e}thodes dans~\cite{gra}.

Dans \cite{mw4}, on a d{\'e}crit pr{\'e}cis{\'e}ment la structure 
multiplicative de la cohomologie {\'e}quivariante des vari{\'e}t{\'e}s de
Bott-Samelson (plus g{\'e}n{\'e}ralement des tours de Bott), puis, en
s'inspirant des m{\'e}thodes utilis{\'e}es par Haibao Duan
dans~\cite{duan}, on a d{\'e}duit 
des r{\'e}sultats de \cite{mw3}  une m{\'e}thode de calcul des
constantes de structure de la cohomologie {\'e}quivariante des vari{\'e}t{\'e}s de
drapeaux (calcul de Schubert {\'e}quivariant). Le but de cet article est
d'effectuer le m{\^e}me travail en $K$-th{\'e}orie. On construit une base de
la $K$-th{\'e}orie {\'e}quivariante des tours de Bott, et on donne une
pr{\'e}sentation par g{\'e}n{\'e}rateurs et relations de ces alg{\`e}bres. On
d{\'e}finit alors r{\'e}cursivement un op{\'e}rateur qui calcule les constantes de 
structure. Gr{\^a}ce aux r{\'e}sultats de \cite{mw3}, on en d{\'e}duit une formule 
pour d{\'e}composer le produit $\hat{\psi}^u\hat{\psi}^v$ sur la base
$\{\hat{\psi}^w\}_{w \in W}$ pour tout couple $(u,v) \in W^2$. 
 
 \bigskip

Les sections~\ref{notations} {\`a} \ref{ktheorie} sont consacr{\'e}es {\`a}
des rappels et des d{\'e}finitions.
Dans la section~\ref{notations}, on fixe les notations sur les
alg{\`e}bres et les groupes de Kac-Moody.
Dans la section~\ref{deftb}, on rappelle la d{\'e}finition et la structure
des tours de Bott et des vari{\'e}tes de Bott-Samelson. Pour plus de
d{\'e}tails, on pourra consulter \cite{livrekumar}, \cite{tours} et \cite{mw4}.
Dans la section~\ref{ktheorie}, on rappelle la d{\'e}finition de la
$K$-th{\'e}orie {\'e}quivariante, la notion de restriction aux points fixes
et la formule de localisation.

\medskip

Dans la section~\ref{sectiontours}, on consid{\`e}re une tour de Bott
$Y$  de dimension complexe $N$ sur laquelle agit le tore compact $D$.
On note $R[D]$ l'anneau des repr{\'e}sentations de $D$.
La $K$-th{\'e}orie $D$-{\'e}quivariante de $Y$, not{\'e}e $K_D(Y)$, est un $R[D]$-module dont 
on construit une base $\{\hat{\mu}_{\epsilon}^D\}_{\epsilon \in \{0,1\}^N}$
(proposition~\ref{propbasektheorie}). On calcule les restrictions
aux points fixes de ces bases (th{\'e}or{\`e}me~\ref{restrictionkb}), et on
donne une pr{\'e}sentation par
g{\'e}n{\'e}rateurs et relations de $K_D(Y)$ (th{\'e}or{\`e}me~\ref{generateursrelationstours}).
Le corollaire~\ref{coror} donne une m{\'e}thode de calcul des constantes
de structure  $r_{\epsilon,
  \epsilon'}^{\epsilon''} \in R[D]$   d{\'e}finies par les relations
$$\hat{\mu}_{\epsilon}^D \hat{\mu}_{\epsilon'}^D
=\sum_{\epsilon'' \in \{0,1\}^N}r_{\epsilon,
  \epsilon'}^{\epsilon''}  \hat{\mu}_{\epsilon''}^D.
$$

\medskip

Dans la section~\ref{sectionBott}, gr{\^a}ce aux r{\'e}sultats
de la section~\ref{sectiontours}, on retrouve la base de la $K$-th{\'e}orie
{\'e}quivariante des vari{\'e}t{\'e}s de Bott-Samelson explicit{\'e}e dans
\cite{mw3} (proposition~\ref{propbasekbs} et
th{\'e}or{\`e}me~\ref{restrictionskbs}), et on donne une m{\'e}thode de
calcul des constantes de structure par rapport {\`a} cette base
(th{\'e}or{\`e}me~\ref{cskbs}).

\medskip

Dans la section~\ref{sectionvd}, on rappelle le lien entre la
$K$-th{\'e}orie des vari{\'e}t{\'e}s de drapeaux et la $K$-th{\'e}orie des
vari{\'e}t{\'e}s de Bott-Samelson gr{\^a}ce au th{\'e}or{\`e}me~\ref{g*k}
d{\'e}montr{\'e} initialement dans \cite{mw3} et dont
on donne une nouvelle d{\'e}monstration. On d{\'e}duit alors des r{\'e}sultats
pr{\'e}c{\'e}dents le th{\'e}or{\`e}me~\ref{cskvd} qui donne
une m{\'e}thode de calcul  des constantes de structure  $q_{u,v}^w \in
R[T]$ d{\'e}finies par les relations
$$\hat{\psi}^u \hat{\psi}^v = \sum_{w \in W}q_{u,v}^w \hat{\psi}^w.
$$

\medskip

La section~\ref{sectionex} est consacr{\'e}e {\`a} quelques
exemples et {\`a} la restriction de nos calculs au cas de la
$K$-th{\'e}orie ordinaire.

\bigskip

\bf Remerciements.
\rm
Je remercie Michel Brion de m'avoir sugg{\'e}r{\'e} d'appliquer mes
r{\'e}sultats au calcul des constantes de structure de la $K$-th{\'e}orie
des vari{\'e}t{\'e}s de drapeaux. Je tiens {\'e}galement {\`a} remercier
Haibao Duan et William Graham pour leurs remarques sur mes
pr{\'e}c{\'e}dents travaux, et Mich{\`e}le Vergne pour sa relecture
attentive de cet article. 

\section{Pr{\'e}liminaires et notations} \label{notations}

\subsection{Alg{\`e}bres de Kac-Moody}

Les d{\'e}finitions et les r{\'e}sultats qui suivent sur les alg{\`e}bres de Kac-Moody
 sont expos{\'e}s dans \cite{km} et \cite{livrekumar}. Soit
 $A=(a_{ij})_{1\leq i,j \leq r}$ une matrice de Cartan
 g{\'e}n{\'e}ralis{\'e}e (c'est-{\`a}-dire 
telle que $a_{ii}=2$, $-a_{ij} \in \mathbb{N}$ si $i \neq j$, et
  $a_{ij}=0$ si et seulement si $a_{ji}=0$). On choisit un triplet $(\mathfrak{h},
  \mathfrak{\pi}, \mathfrak{\pi^{\vee}})$ (unique {\`a} isomorphisme pr{\`e}s), o{\`u}
  $\mathfrak{h}$ est un $\mathbb{C}$-espace vectoriel de dimension $(2r-{ \rm
 rg}(A))$,
  $\mathfrak{\pi} = \{\alpha_{i}\}_{1 \leq i \leq r} \subset \mathfrak{h}^*$, et 
$\mathfrak{\pi^{\vee}} = \{h_i\}_{1 \leq i \leq r} \subset \mathfrak{h}$ sont
  des ensembles d'{\'e}l{\'e}ments lin{\'e}airement ind{\'e}pendants v{\'e}rifiant
  $\alpha_{j}(h_i)=a_{ij}$.  On note aussi $h_i$ par $\alpha^{\vee}_{i}$.
 L'alg{\`e}bre de Kac-Moody
  $\mathfrak{g}=\mathfrak{g}(A)$ est l'alg{\`e}bre de Lie sur $\mathbb{C}$ engendr{\'e}e
  par $\mathfrak{h}$ et par les symboles $e_{i}$ et $f_{i}$ ($1 \leq i \leq r$)
  soumis aux relations $[\mathfrak{h},\mathfrak{h}]=0$, $[h,
  e_{i}]=\alpha_{i}(h)e_{i}$, $[h, f_{i}]=-\alpha_{i}(h)f_{i}$ pour tout $h \in
  \mathfrak{h}$ et tout $1 \leq i \leq r$, $[e_{i}, f_{j}]=\delta_{ij}h_{j}$
  pour tout $1 \leq i,j \leq r$, et  
$$({\rm ad }e_{i})^{1-a_{ij}}(e_{j})=0=({\rm ad
   }f_{i})^{1-a_{ij}}(f_{j})\, ,
   \hspace{0,2 cm} { \rm pour \, tous } \hspace{0,2 cm}  1 \leq i \neq j \leq r.$$

L'alg{\`e}bre $\mathfrak{h}$ s'injecte canoniquement dans $\mathfrak{g}$. On
l'appelle la sous-alg{\`e}bre de Cartan de $\mathfrak{g}$. On a la
d{\'e}composition
$$\mathfrak{g}=\mathfrak{h} \oplus \sum_{\alpha \in
  \Delta_{+}}(\mathfrak{g}_{\alpha} \oplus \mathfrak{g}_{-\alpha}),$$
o{\`u} pour $\lambda \in \mathfrak{h}^*$, $\mathfrak{g}_{\lambda} = \{ x \in
\mathfrak{g} \: {\rm tels \: que }\: [h, x]=\lambda(h)x, \forall h \in
\mathfrak{h} \}$, et o{\`u} on d{\'e}finit $\Delta_{+}$ par $\Delta_{+} = \{ \alpha \in
\sum_{i=1}^{r}\mathbb{N}\alpha_{i} \: {\rm tels \: que } \: \alpha \neq 0
\:{\rm et } \: \mathfrak{g}_{\alpha} \neq 0 \}$. On pose $\Delta=\Delta_{+} \cup
\Delta_{-}$ o{\`u} $\Delta_{-} = -\Delta_{+}$. On appelle $\Delta_{+}$
(respectivement $\Delta_{-}$) l'ensemble des racines positives (respectivement
n{\'e}gatives). Les racines $\{\alpha_{i}\}_{1 \leq i \leq r}$ sont appel{\'e}es les
racines simples. On d{\'e}finit une sous-alg{\`e}bre de Borel $\mathfrak{b}$ de $\mathfrak{g}$ par
$\mathfrak{b}=\mathfrak{h} \oplus \sum_{\alpha \in \Delta_{+}}
\mathfrak{g}_{\alpha}$.

\medskip

Au couple $(\mathfrak{g}, \mathfrak{h})$, on associe le groupe de Weyl 
$W\subset { \rm Aut}(\mathfrak{h}^*)$, engendr{\'e} par les r{\'e}flexions simples 
$\{s_{i}\}_{1  \leq i \leq r}$ d{\'e}finies par 
$$ \forall \lambda \in \mathfrak{h}^*, \,\,
 s_{i}(\lambda)=\lambda-\lambda(h_{i})\alpha_{i} .$$

Si on note $S$ l'ensemble des r{\'e}flexions simples, le couple $(W,S)$ est un
syst{\`e}me de Coxeter. On a donc une notion d'ordre de Bruhat qu'on note 
$u \leq v$ et une notion de longueur
qu'on note $l(w)$. On note $1$ l'{\'e}l{\'e}ment neutre de $W$. Dans le cas fini
(i.e. $W$ fini $\Leftrightarrow
\mathfrak{g}$ de dimension finie),
on note $w_{0}$ le plus grand {\'e}l{\'e}ment de $W$.

On obtient une repr{\'e}sentation de $W$ dans $\mathfrak{h}$ par dualit{\'e}. Plus
pr{\'e}cis{\'e}ment, pour tout $1 \leq i \leq r$, on a : 
$$\forall h \in \mathfrak{h}, \, s_i(h)=h-\alpha_i(h)h_i.$$

Le groupe de Weyl pr{\'e}serve $\Delta$. On pose $R=W\pi$, c'est l'ensemble des
racines r{\'e}elles. On pose $R^{+}=R\cap
\Delta_{+}$, et pour $\beta = w\alpha_{i} \in R^{+}$, on pose
$s_{\beta}=ws_{i}w^{-1} \in W$ (qui est ind{\'e}pendant du choix du couple $(w,
\alpha_{i})$ v{\'e}rifiant $\beta = w\alpha_{i}$) et 
$\beta^{\vee} = wh_{i} \in \mathfrak{h}$.


Pour tout {\'e}l{\'e}ment $w$ de
$W$, on d{\'e}finit l'ensemble $\Delta(w)$ des inversions de $w$ par
 $\Delta(w)=\Delta_{+} \cap w^{-1}\Delta_{-}$.

\medskip

On fixe un r{\'e}seau $\mathfrak{h}_{\mathbb{Z}} \subset \mathfrak{h}$ tel que :

\smallskip

\begin{enumerate}

\item[$(i)$] $\mathfrak{h}_{\mathbb{Z}} \otimes_{\mathbb{Z}}\mathbb{C}=\mathfrak{h}$,

\item[$(ii)$] $h_{i} \in \mathfrak{h}_{\mathbb{Z}}$ pour tout $1 \leq i \leq r$,

\item[$(iii)$] $\mathfrak{h}_{\mathbb{Z}}/ \sum_{i=1}^{r}\mathbb{Z}h_{i}$ est sans
torsion,

\item[$(iv)$] $\alpha_{i} \in \mathfrak{h}_{\mathbb{Z}}^* =
{ \rm Hom}(\mathfrak{h}_{\mathbb{Z}}, \mathbb{Z})$ ($\subset \mathfrak{h}^*$) pour tout 
$1 \leq i \leq r$.
 
\end{enumerate}

\medskip

On choisit des poids fondamentaux $\rho_{i} \in \mathfrak{h}_{\mathbb{Z}}^*$ 
($1 \leq i \leq r$) qui v{\'e}rifient $\rho_{i}(h_{j})=\delta_{i, j}$, pour tout 
$1 \leq i,j \leq r$. On pose $\rho=\sum_{i=1}^{r}\rho_{i}$.

\subsection{Groupes de Kac-Moody et vari{\'e}t{\'e}s de drapeaux}

On note $G=G(A)$ le groupe de Kac-Moody associ{\'e} {\`a} $\mathfrak{g}$ par Kac et
Peterson dans \cite{kp}.  On note $e$ l'{\'e}l{\'e}ment neutre de $G$. 
Dans le cas fini, $G$ est un groupe de Lie semi-simple
complexe connexe et simplement connexe. On note $H \subset B \subset
G$ les sous-groupes de $G$
associ{\'e}s respectivement {\`a} $\mathfrak{h}$ et $\mathfrak{b}$. Soit $K$ la forme
unitaire standard de $G$ et $T=K \cap H$ le tore compact maximal de $K$ associ{\'e} {\`a}
$\mathfrak{h}$. On note $\mathfrak{t} \subset \mathfrak{h}$ 
l'alg{\`e}bre de Lie de $T$. Les racines
$\alpha_i$ et les poids fondamentaux $\rho_i$ appartiennent {\`a}  $i\mathfrak{t}^*$

\medskip

 Soit $N_{G}(H)$ 
le normalisateur de $H$ dans $G$, le groupe quotient
$N_{G}(H)/H$ s'identifie {\`a} $W$. On pose $X=G/B=K/T$. C'est une vari{\'e}t{\'e} de
drapeaux g{\'e}n{\'e}ralis{\'e}e. On fait agir $G$ sur $X$ par
multiplication {\`a} gauche, ce qui induit une action de $B$, $H$ et $T$
sur $X$. L'ensemble des points fixes de $T$ dans $X$ s'identifie {\`a}
$W$.
 Pour $w
\in W$, on d{\'e}finit $C(w)=B \cup BwB$ et pour toute racine simple $\alpha$, on
d{\'e}finit le sous-groupe $P_{\alpha}$ de $G$ par
$P_{\alpha}=C(s_{\alpha})$.
On a la d{\'e}composition de Bruhat $G=\bigsqcup_{w\in W}BwB$ et
si on pose $X_{w}=BwB/B$, $X=\bigsqcup_{w\in W}X_{w}$. Pour tout $w \in W$,
 la cellule de Schubert $X_{w}$ est isomorphe {\`a}
$\mathbb{R}^{2l(w)}$. On obtient ainsi  une d{\'e}composition cellulaire
T-invariante  de $X$ o{\`u} toutes les cellules sont de dimension paire. 

\medskip

Pour tout $w \in W$, la vari{\'e}t{\'e} de Schubert $\overline{X}_{w}$ est l'adh{\'e}rence 
de la cellule $X_{w}$. 
C'est une sous-vari{\'e}t{\'e} irr{\'e}ductible et $T$-invariante de $X$ de dimension r{\'e}elle
$2l(w)$. Les vari{\'e}t{\'e}s de Schubert ne sont pas lisses en g{\'e}n{\'e}ral. Pour tout $w
\in W$, on a la d{\'e}composition suivante : 
$$\overline{X}_{w}=\bigsqcup_{w' \leq w }X_{w'}.$$

\subsection{Le mono{\"\i}de \underline{$W$} }

On d{\'e}finit le mono{\"\i}de $\underline{W}$ comme le mono{\"\i}de engendr{\'e} par les {\'e}l{\'e}ments
$\{\underline{s}_{i}\}_{1 \leq i \leq r}$ soumis aux relations
$\underline{s}_{i}^2=\underline{s}_{i}$ et
aux  relations de tresses  de $W$ :  
$$ \left\{ \begin{array}{cc}
     \underline{s}_{i}^2=\underline{s}_{i} & \\ 
  \underbrace{\underline{s}_{i}\, \underline{s}_{j} \cdots }_{m_{i,j}\, { \rm termes}}=
  \underbrace{\underline{s}_{j}\, \underline{s}_{i} \cdots }_{m_{i,j} \, { \rm
      termes}  }
  & { \rm si }\,
  m_{i,j}<\infty \, ,
   \end{array} \right.$$
o{\`u} $m_{i,j}$ est l'ordre de  $s_is_j$ dans $W$.

\medskip

D'apr{\`e}s l'{\'e}tude g{\'e}n{\'e}rale des alg{\`e}bres de Hecke (voir \cite{hum}), l'ensemble 
$\underline{W}$ s'identifie {\`a} l'ensemble
$W$. Pour un {\'e}l{\'e}ment $w$ de $W$, on note $\underline{w}$ l'{\'e}l{\'e}ment correspondant
dans $\underline{W}$ d{\'e}fini par $\underline{w}=\underline{s}_{i_{1}} \cdots
\underline{s}_{i_{l}}$ si $w=s_{i_{1}} \cdots s_{i_{l}}$ est une d{\'e}composition
r{\'e}duite de $w$, et pour $\underline{v} \in \underline{W}$, on note $v$ l'{\'e}l{\'e}ment
associ{\'e} dans $W$.

\medskip

Dans $\underline{W}$, on a les relations suivantes :

\begin{equation} \label{hecke1}  \left\{ \begin{array}{ll}
      \underline{w}\, \underline{s}_{i}=
\underline{ws_{i}}
 & {\rm si } \, \, ws_{i}>w, \\ 
\underline{w} \, \underline{s}_{i}=\underline{w}
 & {\rm si } \, \, ws_{i}<w.  \end{array} \right.
\end{equation}

\begin{equation} \label{hecke2}
 \left\{ \begin{array}{ll} \underline{s}_{i}\, \underline{w}=\underline{s_{i}w}
 & {\rm si } \, \, s_{i}w>w, \\ 
\underline{s}_{i} \, \underline{w}=\underline{w}
 & {\rm si } \, \, s_{i}w<w.  \end{array} \right.
\end{equation}

\section{Tours de Bott et vari{\'e}t{\'e}s de Bott-Samelson} \label{deftb}

Soit $N \geq 1$ un entier naturel. On pose $\mathcal{E}=\{0,1\}^N$. Pour $\epsilon 
=(\epsilon_1, \epsilon_2, \ldots, \epsilon_N) \in \mathcal{E}$, on
note $\pi_+(\epsilon)$ 
 l'ensemble des entiers $i \in \{1,2,\ldots, N\}$ tels que
$\epsilon_{i}=1$ et $\pi_{-}(\epsilon)$ l'ensemble des entiers $i 
\in \{1,2,\ldots, N\} $
tels que $\epsilon_{i}=0$.
On appelle longueur de $\epsilon$, not{\'e}e
$l(\epsilon)$, le cardinal de $\pi_{+}(\epsilon)$. 
Pour $1\leq i \leq N$, on note $(i) \in \mathcal{E}$
 l'{\'e}l{\'e}ment de $\mathcal{E}$ d{\'e}fini par $(i)_{j}=\delta_{i,j}$.
 On d{\'e}finit  les {\'e}l{\'e}ments  $(\bf{0})$ et $(\bf{1})$ de $\mathcal{E}$
 par $(\bf{0}\rm)_{j}=0$ 
 et $(\bf{1}\rm)_{j}=1$ 
pour tout $j$.  On munit 
$\mathcal{E}$ d'une structure de groupe en
 identifiant $\{0,1\}$ avec $\mathbb{Z}/2\mathbb{Z}$. Pour tout entier $1 \leq n
 \leq N$, on pose $(\overline{n})=(1)+(2)+ \cdots + (n) \in \mathcal{E}$.

 On d{\'e}finit un ordre partiel sur $\mathcal{E}$ par 
$$\epsilon \leq \epsilon' \Leftrightarrow \pi_{+}(\epsilon) \subset
\pi_{+}(\epsilon').$$

\subsection{Tours de Bott} \label{tours}

Les d{\'e}finitions et les r{\'e}sultats des sections~\ref{211} et ~\ref{212}
sont expos{\'e}s plus en d{\'e}tails dans \cite{tours}.

\subsubsection{D{\'e}finition} \label{211}

Les tours de Bott sont des vari{\'e}t{\'e}s complexes compactes et lisses construites de
la mani{\`e}re suivante : 

Soit $\mathbf{L}_2$ un fibr{\'e} en droites holomorphe sur $\mathbb{C}P^1$. On pose
$Y_2=\mathbb{P}(\mathbf{1} \oplus \mathbf{L}_2)$, o{\`u} $\mathbf{1}$ est le fibr{\'e}
en droites 
trivial au dessus de $\mathbb{C}P^1$. La vari{\'e}t{\'e} $Y_2$ est un fibr{\'e} au dessus de
$Y_1=\mathbb{C}P^1$ de fibre $\mathbb{C}P^1$; c'est une surface de
Hirzebruch. On peut it{\'e}rer ce processus {\`a} l'aide de fibr{\'e}s en droites
not{\'e}s $\mathbf{L}_2,  \mathbf{L}_3, \ldots , \mathbf{L}_N$. 
A chaque {\'e}tape, la vari{\'e}t{\'e} $Y_j$ est un
fibr{\'e} au dessus de $Y_{j-1}$ de fibre $\mathbb{C}P^1$. On obtient alors le
diagramme suivant (o{\`u} pour tout $2 \leq j \leq N$, $\mathbf{L}_j$ est un fibr{\'e}
en droites au dessus de $Y_{j-1}$) :  
$$ \begin{array}{clccccc}
 & &  & & & \mathbb{P}(\mathbf{1}\oplus \mathbf{L}_N)=&Y_N \\ 
 &  &  & &  &  \downarrow \pi_N &  \\
  &  & & & &  Y_{N-1}&  \\
 & &   & &\Ddots &  \\
  & & \mathbb{P}(\mathbf{1}\oplus  \mathbf{L}_2)=& Y_2&   & &  \\
  & & \downarrow \pi_2 & &  & &  \\
  & \mathbb{C}P^1 \,\,\, = & \!\! Y_1 & & & &  \\
  &\,\, \downarrow \pi_1  & & &  & &  \\
 \{un \, point\}= &\,\,  Y_0 & & &  & &  
\end{array}$$

\bigskip

A chaque {\'e}tape, on a deux sections particuli{\`e}res $s_j^0 : Y_{j-1} \rightarrow
Y_j$ et $s_j^{\infty} : Y_{j-1} \rightarrow Y_j$ d{\'e}finies par 
$s_j^0(x)=(x,[1,0])$ et $s_j^{\infty}(x)=(x,[0,1])$.

\medskip

Par d{\'e}finition, une tour de Bott de dimension $N$ 
est une famille $\{Y_j, \pi_j, s_j^0,s_j^{\infty}\}_{1 \leq j \leq N}$ issue
d'un diagramme du type pr{\'e}c{\'e}dent.

 On dit que deux tours de Bott 
$\{Y_j, \pi_j, s_j^0,s_j^{\infty}\}_{1 \leq j \leq N}$ et 
$\{{Y'}_j, {\pi'}_j, {s'}_j^0,{s'}_j^{\infty}\}_{1 \leq j \leq N}$
sont isomorphes s'il existe $N$ diff{\'e}omorphismes holomorphes $\{F_j : 
Y_j \rightarrow {Y'}_j \}_{1 \leq j \leq N}$ qui commutent avec les applications 
$\pi_j$, $s_j^0$, $s_j^{\infty}$ et ${\pi'}_j$, ${s'}_j^0$,
${s'}_j^{\infty}$.

\begin{exem}

$\mathbb{C}P^1 \times \cdots \times \mathbb{C}P^1$ (N fois) est une tour de Bott
de dimension $N$.

\end{exem}

\subsubsection{Classes d'isomorphisme des tours de Bott}  \label{212}

On se donne une liste d'entiers $C=\{c_{i,j}\} _{1 \leq
  i < j \leq N}$. On consid{\`e}re $\mathbb{R}^N$ muni de sa base canonique
  $(e_{1}, e_{2}, \ldots , e_{N})$, et on d{\'e}finit $N$ {\'e}l{\'e}ments $v_{1}, v_{2},
  \ldots , v_{N}$ de $\mathbb{R}^{N}$ par les formules suivantes : 
\begin{eqnarray*}
v_{N} & = &-e_{N}, \\
v_{N-1}&=&-e_{N-1}-c_{N-1,N}e_{N},\\
  &\vdots &  \\
v_{1}&=&-e_{1}-c_{1,2}e_{2}- \cdots -c_{1,N}e_{N}.
\end{eqnarray*}

\bigskip

On d{\'e}finit l'{\'e}ventail $\Sigma_{C}$ de $\mathbb{R}^N$ comme la r{\'e}union de tous
les c{\^o}nes engendr{\'e}s par les vecteurs de sous-ensembles $\Lambda$ de $\{ e_{1}, e_{2}, 
\ldots , e_{N}, v_{1}, v_{2}, \ldots , v_{N} \}$ tels que si $e_{i} \in
\Lambda$, alors $v_{i} \notin \Lambda$. On note alors 
$Y_{C}$ la vari{\'e}t{\'e} torique associ{\'e}e {\`a} l'{\'e}ventail $\Sigma_{C}$ (voir 
\cite{cox}, ou \cite{livreaudin} chapitre $6$), c'est-{\`a}-dire le quotient de
  $(\mathbb{C}^2 \setminus {(0,0)})^{N}$ par l'action {\`a} droite de
  $(\mathbb{C}^*)^N$ 
o{\`u} le $i$-{\`e}me
 facteur de $(\mathbb{C}^*)^N$ agit sur  $(\mathbb{C}^2 \setminus {(0,0)})^{N}$ par 
\begin{equation}\label{defaction}
\begin{array}{cc}(z_{1}, w_{1}, \ldots , z_{i-1}, w_{i-1}, z_{i},
  w_{i},z_{i+1}, w_{i+1}, 
 \ldots,  z_{N}, w_{N} )a_{i}=\\
(z_{1}, w_{1}, \ldots , 
z_{i-1}, w_{i-1}, z_{i}a_{i}, w_{i}a_{i},  z_{i+1},
w_{i+1}a_{i}^{c_{i,i+1}},   \ldots ,z_{N}, w_{N}a_{i}^{c_{i,N}}   ).
\end{array}
\end{equation}
 
\medskip

On obtient ainsi une vari{\'e}t{\'e} complexe de dimension $N$. La vari{\'e}t{\'e} $Y_C$ est
 compacte car l'{\'e}ventail
 $\Sigma_{C}$ est complet dans $\mathbb{R}^N$ (i.e. la r{\'e}union des c{\^o}nes de
 $\Sigma_{C}$ est {\'e}gale {\`a} $\mathbb{R}^N$), et lisse car l'{\'e}ventail $\Sigma_{C}$
 est r{\'e}gulier (i.e. les c{\^o}nes de $\Sigma_{C}$ sont engendr{\'e}s par des {\'e}l{\'e}ments du
 r{\'e}seau $\mathbb{Z}^N \subset \mathbb{R}^N$ qui peuvent {\^e}tre compl{\'e}t{\'e}s en une
 base de $\mathbb{Z}^N$). 

On note
$[z_{1}, w_{1}, \ldots ,
z_{N}, w_{N}]$ la classe de $(z_{1}, w_{1}, \ldots , z_{N}, w_{N})$ dans $Y_C$. 

Soit $\epsilon \in \mathcal{E}$. On note $\{i_{1}<
i_{2}<\cdots <i_k \}$ les {\'e}l{\'e}ments de $\pi_+(\epsilon)$. On d{\'e}finit 
 alors une liste d'entiers
$C(\epsilon)=\{d_{i,j}(\epsilon)
\} _{1 \leq
  i < j \leq k}$ par $d_{l,m}(\epsilon)=c_{i_l,i_m}$. En particulier, pour tout entier
$1\leq n \leq N$, on pose $C_n=C((\overline{n}))=\{ c_{i,j} \}_{1 \leq i<j \leq
    n}$.  

Pour tout $2 \leq n \leq N$, 
$Y_{C_n}$ est un fibr{\'e} au dessus de $Y_{C_{n-1}}$ de fibre $\mathbb{C}P^1$. En effet,
on d{\'e}finit un fibr{\'e} en droites  $\mathbf{L}(C_{n-1},c_{1,n},c_{2,n}, \ldots,
c_{n-1,n}) $ sur 
$Y_{C_{n-1}}$ par $\mathbf{L}(C_{n-1},c_{1,n},c_{2,n}, \ldots, c_{n-1,n})=
(\mathbb{C}^2 \setminus {(0,0)})^{n-1}\times_{(\mathbb{C}^*)^{n-1}}\mathbb{C}$, o{\`u}
le $i$-{\`e}me facteur de $(\mathbb{C}^*)^{n-1}$ agit par 
$$((z_{1}, w_{1}, \ldots , z_{n-1}, w_{n-1}),v)a_i=
((z_{1}, w_{1}, \ldots , z_{n-1}, w_{n-1})a_i,a_i^{c_{i,n}}v).$$

Ici l'action de $a_i$ sur $ (z_{1}, w_{1}, \ldots , z_{n-1}, w_{n-1}) $ est donn{\'e}e
par 
$$\begin{array}{cc} (z_{1}, w_{1}, \ldots , z_{i-1}, w_{i-1}, z_{i}, w_{i},z_{i+1}, w_{i+1},
 \ldots,  z_{n-1}, w_{n-1} )a_{i}= \\
(z_{1}, w_{1}, \ldots ,
z_{i-1}, w_{i-1}, z_{i}a_{i}, w_{i}a_{i},  z_{i+1},
w_{i+1}a_{i}^{c_{i,i+1}},  \ldots ,z_{n-1}, w_{n-1}a_{i}^{c_{i,n-1}}
 ). \end{array} $$ 

 \smallskip

On v{\'e}rifie imm{\'e}diatement qu'on a bien $\mathbb{P}(\mathbf{1}\oplus
\mathbf{L}_n)=Y_{C_n}$, o{\`u} $\mathbf{1}$ d{\'e}signe le fibr{\'e} en droites trivial au dessus de 
$Y_{C_{n-1}} $, et o{\`u} $\mathbf{L}_n$ d{\'e}signe le fibr{\'e} 
$\mathbf{L}(C_{n-1},c_{1,n},c_{2,n}, \ldots, c_{n-1,n})$. 
Si on d{\'e}finit $\pi_n : Y_{C_n} \rightarrow Y_{C_{n-1}}$ par 
$$\pi_n([z_{1}, w_{1}, \ldots  ,z_{n-1}, w_{n-1}  , z_{n}, w_{n}  ])
= [z_{1}, w_{1}, \ldots  ,z_{n-1}, w_{n-1} ], $$
la vari{\'e}t{\'e} $Y_{C}$ est alors construite {\`a} l'aide de fibrations
successives de fibres $\mathbb{C}P^1$ selon le diagramme suivant : 
 $$ \begin{array}{clccccc}
 & &  & & & \mathbb{P}(\mathbf{1}\oplus \mathbf{L}_N)=&Y_{C_N} = Y_C \\ 
 &  &  & &  &  \downarrow \pi_N &  \\
  &  & & & &  Y_{C_{N-1}}&  \\
 & &   & &\Ddots &  \\
  & & \mathbb{P}(\mathbf{1}\oplus  \mathbf{L}_2)=& Y_{C_2}&   & &  \\
  & & \downarrow \pi_2 & &  & &  \\
  & \mathbb{C}P^1 \,\,\, = & \!\! Y_{C_1} & & & &  \\
  &\,\, \downarrow \pi_1  & & &  & &  \\
 \{un \, point\}= &\,\,  Y_0 & & &  & &  
\end{array}$$

\bigskip

Si on d{\'e}finit $s_j^0 : Y_{j-1} \rightarrow
Y_j$ et $s_j^{\infty} : Y_{j-1} \rightarrow Y_j$ comme dans la section
pr{\'e}c{\'e}dente, alors la famille $\{Y_{C_j}, \pi_j, s_j^0, s_j^{\infty}  \}_{1 \leq
  j \leq N}$
est donc une tour de Bott de dimension $N$, et on obtient ainsi une application : 
\begin{equation} \label{isomtours}
\mathbb{Z}^{N(N-1)/2} \rightarrow \{ \mbox{ classes  d'isomorphisme 
de  tours  de   Bott 
   de  dimension } \, N \},
\end{equation}
qui {\`a} $C \in \mathbb{Z}^{N(N-1)/2}$ associe $Y_C$,
et pour toute tour de Bott $Y_C$ dans l'image de~\ref{isomtours}, une application : 
\begin{equation} \label{isomfibr{\'e}s}
\mathbb{Z}^{N} \rightarrow \{ \mbox{ classes  d'isomorphisme 
de  fibr{\'e}s  en   droites 
   holomorphes  sur }\,  Y_C \},
\end{equation}
qui {\`a} $(m_1,m_2, \ldots, m_N) \in \mathbb{Z}^{N}$ associe le fibr{\'e} 
$\mathbf{L}(C,m_1,m_2, \ldots, m_N)$ sur $Y_C$.

\smallskip

Le r{\'e}sultat suivant est alors prouv{\'e} dans \cite{tours} : 

\begin{prop}

Les applications~\ref{isomtours} et \ref{isomfibr{\'e}s} sont des bijections.

\end{prop}

\subsubsection{Actions de groupes sur les tours de Bott}

Pour tout $1 \leq i \leq N$, on fait agir
$D_{\mathbb{C}}=(\mathbb{C}^*)^N$ (d'alg{\`e}bre de Lie 
$\mathfrak{d}_{\mathbb{C}} \simeq \mathbb{C}^N$) 
sur $Y_{C_i}$ par 
\begin{equation} \label{actiontours}
\begin{array}{c}
(e^{\lambda_{1}(d)},e^{\lambda_{2}(d)}, \ldots , e^{\lambda_{N}(d)})
[z_{1}, w_{1}, z_{2}, w_{2}, \ldots ,z_{i},
w_{i}]=\\
{[}z_{1}, e^{-\lambda_{1}(d)} w_{1},z_{2}, e^{-\lambda_{2}(d)} w_{2},
 \ldots ,z_{i}, e^{-\lambda_{i}(d)} w_{i}{]},
\end{array}
\end{equation}
o{\`u} $\lambda_k \in \mathfrak{d}_{\mathbb{C}}^*$ est d{\'e}finie par 
$\lambda_k((d_1,d_2, \ldots , d_N))=d_k$. 

\smallskip



 L'action de $D_{\mathbb{C}}$ sur $Y_{C_i}$ 
induit une action de $D=(S^1)^N$ (d'alg{\`e}bre de Lie $\mathfrak{d} \simeq
i\mathbb{R}^N \subset\mathfrak{d}_{\mathbb{C}} $) sur $Y_{C_i}$. 
Les $\lambda_k$ sont dans $i\mathfrak{d}^*$.

\bigskip

Pour tout $2 \leq i \leq N$, on fait agir $D_{\mathbb{C}}$ (et par
restriction $D$) sur
$\mathbf{L}_i$ par  
\begin{equation} \label{actionfibres}
\begin{array}{c}(e^{\lambda_1(d)},e^{\lambda_2(d)}, \ldots , 
   e^{\lambda_{N}(d)})[z_1, w_1, z_2, w_2, \ldots ,z_{i-1},
w_{i-1},v]= \\ 
 {[}z_1, e^{-\lambda_1(d)} w_{1},z_2, e^{-\lambda_2(d)} w_2,
 \ldots ,z_{i-1}, e^{-\lambda_{i-1}(d)}
 w_{i-1},e^{-\lambda_i(d)}v{]}.
\end{array}
\end{equation}

De plus, on note $\mathbf{L}_1 \simeq \mathbb{C}$ le fibr{\'e} en droites
trivial sur le point, et on fait agir $D_{\mathbb{C}}$ sur
  $\mathbf{L}_1$ par  
$$(e^{\lambda_1(d)},e^{\lambda_2(d)}, \ldots ,
e^{\lambda_N(d)})v=e^{-\lambda_1(d)}v.
$$

\medskip

Le lemme suivant est imm{\'e}diat : 

\begin{lemm} \label{isomequivariant}
Pour tout $1 \leq i \leq N$, l'action de $D_{\mathbb{C}}$ sur
$Y_{C_i}=\mathbb{P}(\mathbf{1}\oplus  \mathbf{L}_i)$ est induite par
l'action de $D_{\mathbb{C}}$ sur $\mathbf{L}_i$. 

De plus, les applications $\pi_i$ sont $D_{\mathbb{C}}$-{\'e}quivariantes.
\end{lemm}

\medskip
 
Dans toute la suite on suppose donn{\'e}e une liste d'entiers $C$ et on note $Y$ 
(au lieu de $Y_{C}$) la tour de Bott associ{\'e}e {\`a} $C$.

\subsubsection{D{\'e}composition cellulaire} On d{\'e}finit une d{\'e}composition cellulaire de
$Y$ index{\'e}e par $\mathcal{E}$ de la mani{\`e}re suivante : 

 Pour $\epsilon \in \mathcal{E}$, on note $Y_{\epsilon}
\subset Y$ l'ensemble des classes $ [z_{1}, w_{1}, \ldots ,z_{N}, w_{N}] $
 qui v{\'e}rifient pour tout entier $i$ compris entre
$1$ et $N$ 
$$ w_{i}=0
  { \rm\,\,\, si} \hspace{0,15 cm} \epsilon_{i} =0 \, , \,\,\,\,\,\,\,
 w_{i}\neq 0
  {\,\,\, \rm si }\hspace{0,15 cm} \epsilon_{i} =1.$$

\smallskip

On v{\'e}rifie imm{\'e}diatement que cette d{\'e}finition est bien compatible avec
l'action de $(\mathbb{C}^*)^{N}$ sur
$(\mathbb{C}^2\setminus{(0,0)})^N$
d{\'e}finie par l'{\'e}quation~\ref{defaction}. 

\medskip

 Pour $\epsilon \in \mathcal{E}$ et $1 \leq k<l \leq N$, on pose 
$$\displaystyle{ c_{k,l} (\epsilon) =-c_{k,l}+ \!\!\!\!\!\!
\sum_{ {\tiny \begin{array}{cc} {k<i_1< \cdots < i_m<l}  \\ 
 m>0, i_{j} \in \pi_{+}(\epsilon) \end{array}}}}\!\!\!\!\!\!\!\!
\!\!\!\!\!\!\!\!(-1)^{m+1} c_{k,i_1}c_{i_1,i_2}
\cdots c_{i_m,l}.$$

Ces nombres sont {\'e}galement donn{\'e}s par r{\'e}currence sur $l-k > 0$ par les
relations suivantes : 
\begin{equation} \label{c}
\left\{ \begin{array}{ll} c_{k,k+1} (\epsilon)=-c_{k,k+1}
 & { \rm pour} \hspace{0,15 cm} 1 \leq k \leq N-1, \\ 
c_{k,l}(\epsilon) = -c_{k,l}- 
\sum_{ {\tiny \begin{array}{cc} {k<m<l}  \\ 
 m \in \pi_{+}(\epsilon) \end{array}}}c_{m,l}c_{k,m}(\epsilon)
 & { \rm pour }\hspace{0,15 cm}1 \leq k <l \leq N,
\end{array}\right. \end{equation}
o{\`u}, par convention, $\sum_{\emptyset}=0$.

\medskip

 Pour $\epsilon \in \mathcal{E}$ et $i \in \{ 1, 2, \ldots , N \}$, on d{\'e}finit
 $\lambda_{i}(\epsilon) \in \mathfrak{d}_{\mathbb{C}}^* $ par
$$\lambda_{i}(\epsilon)=(-1)^{\epsilon_{i}+1} \big( \lambda_{i}+\sum_{j<i , j \in
  \pi_{+}(\epsilon)}\!\!\!\!\!\!c_{j,i}(\epsilon)\lambda_{j} \big).$$

On d{\'e}montre alors facilement la proposition suivante : 

\begin{prop}  \label{decompositiontours}
    
\indent

\begin{enumerate}

\item[$(i)$] Pour tout $\epsilon \in \mathcal{E}$, $Y_{\epsilon}$ est un
espace affine complexe de dimension $l(\epsilon)$ stable sous l'action 
lin{\'e}aire du tore $D_{\mathbb{C}}$. 

\item[$(ii)$] Pour tout $\epsilon \in \mathcal{E}$, $\overline{Y}_{\epsilon} = 
\coprod_{\epsilon' \leq \epsilon} Y_{\epsilon'}$.

\item[$(iii)$] $Y = \coprod_{\epsilon \in
  \mathcal{E}} Y_{\epsilon}$.

\item[$(iv)$] Pour tout $\epsilon  \in \mathcal{E}$, la sous-vari{\'e}t{\'e}
$\overline{Y}_{\epsilon}$ s'identifie {\`a} la vari{\'e}t{\'e} $Y_{C(\epsilon) }$ 
 et est donc une sous-vari{\'e}t{\'e} irr{\'e}ductible lisse de $Y$.

\end{enumerate}

\end{prop}

De plus, nous allons avoir besoin du lemme suivant dont la
d{\'e}monstration est imm{\'e}diate :

\begin{lemm} \label{pointsfixestours}

\indent

\begin{enumerate}
    
\item[$(i)$] L'ensemble $Y^{D}$ des points fixes de $Y$ sous l'action de
$D$ est constitu{\'e} des $2^N$ points :
$$[z_{1}, w_{1},  z_{2}, w_{2}, \ldots, z_{N}, w_{N}], 
\hspace{0,15 cm} o\grave{u} \hspace{0,15 cm}
(z_{i},w_{i}) \in \{ (1,0), (0,1) \}.$$

On identifie donc $Y^{D}$ avec $\mathcal{E}$ en identifiant $(1,0)$
avec $0$ et $(0,1)$ avec $1$. Le point fixe $\epsilon \in Y^{D}$ est
l'unique point fixe de $Y_{\epsilon}$.

\item[$(ii)$] Soit $(\epsilon,\epsilon') \in \mathcal{E}^2$, alors : 
$$ \epsilon \in \overline{Y}_{\epsilon'} \Leftrightarrow \epsilon \leq \epsilon',$$
et dans ce cas si on note $T_{\epsilon'}^{\epsilon}$ l'espace tangent {\`a}
$\overline{Y}_{\epsilon'}$ en $\epsilon$, les poids de la repr{\'e}sentation de
$D$ dans $T_{\epsilon'}^{\epsilon}$ induite par l'action de
$D$ sur
$Y$ sont les $\{ \lambda_{i}(\epsilon) \}_{i \in \pi_{+}(\epsilon') }$. 

\end{enumerate}

\end{lemm}

\subsection{Vari{\'e}t{\'e}s de Bott-Samelson} \label{BS}

On utilise les notations de la section~\ref{notations}. Pour plus de
d{\'e}tails sur les vari{\'e}t{\'e}s de Bott-Salmelson, on pourra consulter
\cite{livrekumar}.

\subsubsection{D{\'e}finition} \label{def}

Consid{\'e}rons une suite de $N$ racines
simples $\mu_{1}$, \ldots , $\mu_{N}$ non n{\'e}cessairement distinctes. 
On d{\'e}finit 
$$\Gamma(\mu_{1}, \ldots ,\mu_{N})=P_{\mu_1} \times_{B} 
P_{\mu_2} \times_{B} \cdots \times_{B} P_{\mu_N}/B,$$
comme l'espace des orbites de $B^N$ dans $P_{\mu_1} \times P_{\mu_2} 
\times \cdots \times P_{\mu_N}$,
sous l'action {\`a} droite de $B^N$ d{\'e}finie par 
$$(g_{1}, g_{2}, \ldots , g_{N})
(b_{1}, b_{2}, \ldots , b_{N}) = 
(g_{1}b_{1},b_{1}^{-1} g_{2}b_{2}, \ldots ,b_{N-1}^{-1} g_{N}b_{N}),\,
 b_{i} \in B, \, g_{i} \in P_{\mu_i}.$$

\medskip

On obtient ainsi  une vari{\'e}t{\'e}
projective irr{\'e}ductible et lisse. On note
$[g_{1}, g_{2}, \ldots , g_{N}]$ la classe de $(g_{1}, g_{2},
\ldots , g_{N})$ dans $\Gamma(\mu_{1}, \ldots ,\mu_{N})$. On note $g_{\mu_{i}}
\in P_{\mu_i}$ un repr{\'e}sentant quelconque de la r{\'e}flexion de 
$N_{P_{\mu_i}}(H)/H \simeq \mathbb{Z}/2\mathbb{Z}$.

\bigskip

On d{\'e}finit une action {\`a} gauche de $B$ sur $\Gamma$ par  
$$b[g_{1}, g_{2}, \ldots ,g_{N}]=[bg_{1}, g_{2}, \ldots ,g_{N}],\hspace{0,1 cm} b \in B,
\hspace{0,1 cm}  g_{i} \in P_{\mu_i}.$$

Par restriction, on obtient ainsi une action de $H$ et de $T$.

\subsubsection{D{\'e}composition cellulaire}

 Pour $\epsilon \in \mathcal{E}$, on note $\Gamma_{\epsilon}
\subset \Gamma$ l'ensemble des classes $[g_{1},
g_{2}, \ldots , g_{N}]$ qui v{\'e}rifient pour tout entier $i$ compris entre
$1$ et $N$ 
$$ g_{i} \in B
  {\,\,\, \rm si} \hspace{0,15 cm} \epsilon_{i} =0 \, , \,\,\,\,\,\,\,\,
 g_{i} \notin B
  {\,\,\, \rm si }\hspace{0,15 cm} \epsilon_{i} =1.
$$

On v{\'e}rifie imm{\'e}diatement que cette d{\'e}finition est bien compatible avec
l'action de $B^N$.

\medskip
 
 Pour $\epsilon \in \mathcal{E}$ et $1 \leq i\leq N$, on d{\'e}finit 
$$\displaystyle{v_{i}(\epsilon) =\prod_{ {\tiny \begin{array}{cc}  1\leq k \leq i, \\ 
 k \in \pi_{ +}(\epsilon) \end{array}} 
  }\!\!\!\!\!\!s_{\mu_{k}}},$$
o{\`u}, par convention, $\prod_{\emptyset}=1$. On pose $v(\epsilon)=v_{N}(\epsilon)$. 
 Ce sont des {\'e}l{\'e}ments de $W$.




\medskip

De plus, on pose 
$\alpha_{i}(\epsilon)=v_{i}(\epsilon)\mu_{i}$, c'est une racine.

\medskip

On d{\'e}finit de m{\^e}me $$\displaystyle{\underline{v}(\epsilon) = 
\prod_{ {\tiny \begin{array}{cc}  1\leq k \leq N, \\ 
 k \in \pi_{ +}(\epsilon) \end{array}} } 
 \!\!\!\!\!\!  \underline{s_{\mu_{k}}} \in \underline{W}}.$$ 

\bigskip
\medskip

On d{\'e}montre alors facilement la proposition suivante : 

\begin{prop} \label{structureBS}
    
\indent

\begin{enumerate}

\item[$(i)$] Pour tout $\epsilon \in \mathcal{E}$, $\Gamma_{\epsilon}$ est un
espace affine complexe de dimension $l(\epsilon)$ stable sous l'action de $B$,
et cette action induit une action lin{\'e}aire du tore $H$ sur $\Gamma_{\epsilon}$. 

\item[$(ii)$] Pour tout $\epsilon \in \mathcal{E}$, $\overline{\Gamma}_{\epsilon} = 
\coprod_{\epsilon' \leq \epsilon} \Gamma_{\epsilon'}$.

\item[$(iii)$] $\Gamma = \coprod_{\epsilon \in
  \mathcal{E}} \Gamma_{\epsilon}$.

\item[$(iv)$] Pour tout $\epsilon  \in \mathcal{E}$,
  $\overline{\Gamma}_{\epsilon}$ s'identifie 
{\`a} la vari{\'e}t{\'e} $\Gamma(\mu_{i}, i \in \pi_{+}(\epsilon))$ et est donc une
sous-vari{\'e}t{\'e} irr{\'e}ductible lisse de $\Gamma$.

\end{enumerate}   

\end{prop}

Soit $\Gamma^T$ l'ensemble des points fixes de $T$ dans $\Gamma$, on peut
identifier $\Gamma^T$ avec $\mathcal{E}$ gr{\^a}ce au lemme suivant :

\begin{lemm} \label{pointsfixesBS}

\indent 
    
\begin{enumerate}

\item[$(i)$] L'ensemble $\Gamma^T$ est
constitu{\'e} des $2^N$ points :
$$[g_{1}, g_{2}, \ldots, g_{N}], \,\, \mbox{o{\`u}} \,\,\,
g_{i} \in \{ e, g_{\mu_{i}} \}.$$

On identifie donc $\Gamma^T$ avec $\mathcal{E}$ en identifiant $e$
avec $0$ et $g_{\mu_{i}}$ avec $1$.

\item[$(ii)$] Pour tout $\epsilon \in \mathcal{E}$, le point fixe $\epsilon$ est l'unique
point fixe de $\Gamma_{\epsilon}$.

\end{enumerate}

\end{lemm}

\subsubsection{Les vari{\'e}t{\'e}s de Bott-Samelson sont des tours de Bott}
\label{toursbott}

Soit $\tau : \mathfrak{d}_{\mathbb{C}}^* \rightarrow \mathfrak{h}^*$
l'application  d{\'e}finie par
$\tau(\lambda_{i}) = \mu_{i}$. Elle envoie $\mathfrak{d}^*$ dans $\mathfrak{t}^*$.
Soit $\mathfrak{s}^* \subset \frak{t}^*$ l'image de  $\mathfrak{d}^*$ par 
$\tau$ et soit 
$\mathfrak{s}_{\mathbb{C}}^* \subset \frak{h}^*$ 
l'image de $\mathfrak{d}_{\mathbb{C}}^*$. On a les deux suites suivantes (o{\`u} les
premi{\`e}res fl{\`e}ches sont surjectives et les deuxi{\`e}mes injectives) :  
 $$\xymatrix{ \mathfrak{d}^* \ar@{->>}[r]^{\tau }\,\,\, & 
\, \frak{s}^*\,\,
 \ar@{^{(}->}[r]&   \mathfrak{t}^* }, $$
$$\xymatrix{ \mathfrak{d}_{\mathbb{C}}^* \ar@{->>}[r]^{\tau }\,\,\, & 
\, \mathfrak{s}_{\mathbb{C}}^*\,\,
 \ar@{^{(}->}[r]&   \mathfrak{h}^* }.
$$

On en d{\'e}duit la suite suivante sur les tores complexes : 
$$\xymatrix{ H \ar@{->>}[r]^{ }\,\,\, & 
\, S_{\mathbb{C}}\,\,
 \ar@{^{(}->}[r]^{\gamma }&   D_{\mathbb{C}} }. $$

On continue {\`a} noter $\gamma$ le morphisme de $H$ dans $D_{\mathbb{C}}$ ainsi
d{\'e}fini.
On a {\'e}galement la suite suivante sur les tores compacts : 
$$\xymatrix{ T \ar@{->>}[r]^{ }\,\,\, & 
\, S\,\,
 \ar@{^{(}->}[r]^{\gamma }&   D }. $$

 L'action de $h \in H$ sur $\Gamma$ {\'e}tant donn{\'e}e par la
formule 
$$h[g_1,g_2, \ldots, g_N]=[hg_1h^{-1},hg_2h^{-1}, \ldots, ,hg_Nh^{-1}],$$
le tore $H$ agit sur $\Gamma$ via son image $S_{\mathbb{C}}$.
 
\bigskip

Si $N=1$, $\Gamma(\mu_1)$ est isomorphe {\`a} $\mathbb{C}P^1$. De plus, si
$N \geq 2$, $\Gamma(\mu_1, \ldots, \mu_N)$ est une fibration au dessus
de $\Gamma(\mu_1, \ldots, \mu_{N-1})$ de fibre $\mathbb{C}P^1$. En
fait, on peut identifier la vari{\'e}t{\'e} de Bott-Samelson $\Gamma$ {\`a} une tour de
Bott gr{\^a}ce {\`a} la proposition suivante (voir \cite{tours} et \cite{mw4}
pour plus de d{\'e}tails) : 

\begin{prop}

Il existe un diff{\'e}omorphisme (de vari{\'e}t{\'e} $C^{\infty}$) $\phi$ entre 
la vari{\'e}t{\'e}  $\Gamma
(\mu_1,\mu_2, \ldots ,\mu_N)$ et la
tour de Bott $Y_C$, avec $C=\{c_{j,k} \}_{1 \leq j < k \leq N}$, o{\`u} $c_{j,k}$ est d{\'e}fini par 
$c_{j,k}=\mu_{k}(\mu^{\vee}_j)$.

La liste $C$ est donc un ensemble de nombres de Cartan d{\'e}pendant du choix de la
suite $\mu_1,\mu_2, \ldots ,\mu_N$.

 De plus, 
$$\forall (h,x) \in H \times \Gamma, \,
\phi(hx)= \gamma(h)\phi(x),  $$
et l'action de $H$ sur $\Gamma$ s'identifie donc {\`a} celle d'un sous-tore de 
$D_{\mathbb{C}}$ sur $Y_{C}$.

\end{prop}

 On v{\'e}rifie, gr{\^a}ce {\`a} la construction 
de $\phi$ (voir \cite{mw4}), que $\phi$ 
 envoie $\Gamma_{\epsilon}$ sur $Y_{\epsilon}$ et le point $\epsilon \in \Gamma^T$ sur 
le point $\epsilon \in Y^{D}$.

Le tore $T$ agit sur 
$\Gamma$ via son image $S$, et on a 
$$\forall (t,x) \in T \times \Gamma, \,
\phi(tx)= \gamma(t)\phi(x).  $$

L'action de $T$ sur $\Gamma$ s'identifie donc {\`a} celle d'un sous-tore de 
$D$ sur $Y_{C}$.




\section{K-th{\'e}orie {\'e}quivariante} \label{ktheorie}

\subsection{D{\'e}finitions}

 Soit $U_{\mathbb{C}}$ un tore complexe 
 d'alg{\`e}bre de Lie $\mathfrak{u}_{\mathbb{C}}$,
et soit $U   \subset U_{\mathbb{C}}$ le tore compact maximal de
 $U_{\mathbb{C}}$.
  On note 
 $\mathfrak{u} \subset \mathfrak{u}_{\mathbb{C}} $ l'alg{\`e}bre de Lie de $U$.
 On note $X[U]$ 
le groupe des caract{\`e}res de $U$, et on pose $R[U]=\mathbb{Z}[X[U]]$. On note
 $Q[U]$ le corps des fractions de $R[U]$.

Pour tout poids entier $\alpha \in i\mathfrak{u}^* \subset 
\mathfrak{u}_{\mathbb{C}}^*$, on note $e^{\alpha} : U \rightarrow S^1$ le
caract{\`e}re correspondant.

\bigskip

Soit $Z$ un espace topologique compact muni d'une action continue de $U$. 
On d{\'e}finit la $K$-th{\'e}orie $U$-{\'e}quivariante de $Z$ comme le groupe construit {\`a}
partir du semi-groupe des classes
d'isomorphisme de fibr{\'e}s vectoriels complexes de dimension finie
$U$-{\'e}quivariants au 
dessus de $Z$. On munit ce groupe
d'une structure d'anneau d{\'e}finie {\`a} l'aide du produit tensoriel. De
plus, comme la  $K$-th{\'e}orie $U$-{\'e}quivariante du point s'identifie {\`a} $R[U]$, on
obtient une structure de 
$R[U]$-alg{\`e}bre qu'on note $K_U(Z)$.

\medskip

Toute application $g : Z_1 \rightarrow Z_2$ continue et $U$-{\'e}quivariante d{\'e}finit
un morphisme de $R[U]$-alg{\`e}bre $g^* : K_{U}(Z_2) \rightarrow K_{U}(Z_1)$. En
particulier, l'inclusion $Z^U \subset Z$ d{\'e}finit un morphisme 
$i_U^*  : K_{U}(Z) \rightarrow K_{U}(Z^U)$ appel{\'e} restriction aux points fixes. 
 Si l'ensemble des points fixes $Z^U$ est discret, $  K_{U}(Z^U) $ 
s'identifie de mani{\`e}re {\'e}vidente {\`a} $F(Z^U; R[U])$ la
$R[U]$-alg{\`e}bre
des fonctions de $Z^U$ {\`a} valeurs dans $R[U]$ munie de l'addition et de la
multiplication point par point. On obtient alors un morphisme $i_U^* : 
  K_{U}(Z)  \rightarrow  F(Z^U; R[U]) $.

 On note~$*$ l'involution de
$K_U(Z)$ d{\'e}finie par la dualit{\'e} des fibr{\'e}s, et on note de la m{\^e}me fa{\c c}on
l'involution de $R[U]$ d{\'e}finie sur les caract{\`e}res par
$*(e^{\alpha})=e^{-\alpha}$, ce qui induit une involution de $F(Z^U ; R[U])$.
 Pour tout {\'e}l{\'e}ment $\tau \in K_{U}(Z)$,
$$*i_U^*(\tau)=i_{U}^*(*\tau).$$

\subsection{Formule de localisation}

On suppose que $Z$ est une vari{\'e}t{\'e} complexe projective lisse
de dimension 
$n$ munie d'une action de 
$U_{\mathbb{C}}$. Cette action induit alors une action de $U$ sur $Z$. 

La vari{\'e}t{\'e} $Z$ {\'e}tant lisse, le groupe de Grothendieck
$K_{0}(U_{\mathbb{C}},Z)$ construit {\`a} partir des classes 
d'isomorphisme de faisceaux 
$U_{\mathbb{C}}$-{\'e}quivariants coh{\'e}rents sur $Z$ est isomorphe au
groupe de Grothendieck 
$K^{0}(U_{\mathbb{C}}, Z)$ construit {\`a} partir des classes 
d'isomorphisme de faisceaux $U_{\mathbb{C}}$-{\'e}quivariants localement libres sur
$Z$ (voir \cite{ginzburg}, chapitre $5$). On identifie donc ces deux groupes.

On a un morphisme canonique : $K^{0}(U_{\mathbb{C}}, Z) \rightarrow K_U(Z)$. On
suppose que ce morphisme est un isomorphisme (c'est le cas pour les tours de Bott
et les vari{\'e}t{\'e}s de drapeaux dans le cas fini), et on identifie ces deux groupes.

 Pour toute sous-vari{\'e}t{\'e} $U_{\mathbb{C}}$-invariante  $Z'$ et tout 
 faisceau $\mathcal{F} \in  K^{0}(U_{\mathbb{C}}, Z) $, l'action de
 $U_{\mathbb{C}}$ sur $Z$ induit une action de $U_{\mathbb{C}}$ sur 
${\rm H}^k(Z',\mathcal{F}_{/Z'})$, et on d{\'e}finit 
$\chi(Z', \mathcal{F}) \in R[U]$ par 
$$\forall u \in U, \hspace{0,1 cm} \chi(Z',
\mathcal{F})(u)=\sum_{k}(-1)^k{ \rm Tr }
 (u;  { \rm H}^k(Z',\mathcal{F}_{/Z'})).$$

On suppose de plus que $Z^U$ est fini. 
En chaque point fixe $m \in Z^U$, on note
$(\alpha_1^m, \ldots ,\alpha_n^m) \in (i\mathfrak{u}^*)^n 
 \subset (\mathfrak{u}_{\mathbb{C}}^*)^n$ les poids de la
repr{\'e}sentation de $U$ dans $T_mZ$, l'espace tangent {\`a} $Z$ en $m$. 
Dans ce cas, la formule $5.11.9$ de \cite{ginzburg} 
s'{\'e}crit de la mani{\`e}re suivante : 

\begin{prop} \label{pointsfixesab}
 
Pour tout faisceau $\mathcal{F}$ localement libre et
$U_{\mathbb{C}}$-{\'e}quivariant au dessus de $Z$, $\chi(Z,\mathcal{F})  $ se
calcule gr{\^a}ce {\`a} la formule suivante  : 

$$\chi(Z,\mathcal{F})=\sum_{m \in Z^U}\frac{i_{U}^*(\mathcal{F})(m) }
{\prod_{1 \leq i \leq n}(1-e^{-\alpha_i^m}) }.$$

\end{prop}

\section{K-th{\'e}orie {\'e}quivariante des tours de Bott} \label{sectiontours}

On reprend les notations de la section~\ref{tours}. Soit $N\geq 1$ un entier
naturel, et soit $C=\{c_{i,j}\}_{1 \leq i<j \leq N}$ une liste
d'entiers. On pose $Y=Y_{C}$.

On montre par r{\'e}currence sur la dimension de $Y$ (voir \cite{kkk} o{\`u} le 
r{\'e}sultat est d{\'e}montr{\'e} dans le cas particulier des vari{\'e}t{\'e}s de 
Bott-Samelson) que le morphisme canonique  $K^{0}(D_{\mathbb{C}}, Y)
\rightarrow K_D(Y)$ est un isomorphisme. 
Dans la suite, on identifie donc ces deux groupes.

\bigskip

Dans \cite{mw3}, on a calcul{\'e} les restrictions aux points fixes
de la base de la $K$-th{\'e}orie {\'e}quivariante des vari{\'e}t{\'e}s de
Bott-Samelson d{\'e}finie par dualit{\'e} par rapport {\`a} la
d{\'e}composition cellulaire construite dans la section~\ref{BS}. 
Ici, on construit d'abord une base de $K_D(Y)$ {\`a} l'aide de
fibr{\'e}s en droites (formule~\ref{muproduit}), on calcule les
restrictions aux points fixes de cette base
(th{\'e}or{\`e}me~\ref{restrictionkb}),  et on en d{\'e}duit que c'est
la base duale par rapport {\`a} la d{\'e}composition cellulaire construite
dans la section~\ref{tours} (th{\'e}or{\`e}me~\ref{ktheorietours}). 
On retrouvera alors dans la section~\ref{sectionBott} la base de la
$K$-th{\'e}orie {\'e}quivariante des vari{\'e}t{\'e}s de Bott-Samelson utilis{\'e}e dans
\cite{mw3}. On va de plus expliciter la structure multiplicative de
ces alg{\`e}bres.

\subsection{Construction d'une base}

La structure de la K-th{\'e}orie $D$-{\'e}quivariante de $Y$ est donn{\'e}e par la
proposition suivante :  

\begin{prop} \label{propbasektheorie}

\indent

\begin{enumerate}

\item[$(i)$] La  $K$-th{\'e}orie $D$-{\'e}quivariante de $Y^D$ s'identifie {\`a}
$F(\mathcal{E};R[D])$.

\item[$(ii)$] La restriction aux points fixes $i_{D}^*$ : $K_D(Y)
 \rightarrow F(\mathcal{E};R[D])$ est injective.

\item[$(iii)$] La $K$-th{\'e}orie $D$-{\'e}quivariante de $Y$ est un $R[D]$-module
  libre de rang $2^N$.
\end{enumerate}

\end{prop}

\begin{proof}

Le point $(i)$ est imm{\'e}diat.

\smallskip

Le point $(ii)$ est une cons{\'e}quence de $(iii)$. En effet, d'apr{\`e}s le th{\'e}or{\`e}me de
localisation (voir \cite{kequivariante}), le morphisme : 
$K_D(Y)\otimes_{R[D]}Q[D] \rightarrow F(\mathcal{E}; Q[D])$ induit par
$i_D^*$ est un isomorphisme. De plus, comme $K_D(Y)$ est un $R[D]$-module
libre, $K_{D}(Y)$ s'injecte dans $  K_{D}(Y)\otimes_{R[D]}Q[D] $, et on a donc
le diagramme commutatif suivant qui prouve que $i_{D}^*$ est injective : 
$$\xymatrix{
     K_{D}(Y)\ar[dd]^{i_D^*} \ar@{^{(}->}[rr] & &
    K_{D}(Y)\otimes_{R[D]}Q[D]\ar[dd]^{\simeq} \\
     \\
     F(\mathcal{E};R[D]) \ar[rr] & & 
   F(\mathcal{E};Q[D]) }
  \\
$$

\bigskip

Pour d{\'e}montrer $(iii)$, on va expliciter une base 
$\{ \hat{\mu}_{\epsilon}^D \}_{\epsilon \in \mathcal{E}}$
du $R[D]$-module $K_D(Y)$. On proc{\`e}de par r{\'e}currence sur $N \geq 1$.

\smallskip

 Pour $N=1$,
$K_{S^1}(\mathbb{P}^1)$ est un $R[S^1]$-module libre engendr{\'e} par le fibr{\'e} en
droites 
trivial $\mathbf{1}$ et par le fibr{\'e} $\mathbf{E}$
d{\'e}fini comme le fibr{\'e} en droites tautologique sur  $\mathbb{P}^1$,
soumis {\`a} l'unique relation 
$$\mathbf{E}^2=(1+e^{-\lambda_1})\mathbf{E}-e^{-\lambda_1}\mathbf{1},$$
 (voir
\cite{atiyah}, corollaire 2.2.2). On pose $ \hat{\mu}_{(0)}^D= \mathbf{E}$ et 
$ \hat{\mu}_{(1)}^D= \mathbf{1} - \mathbf{E}$.

\smallskip

On suppose le r{\'e}sultat v{\'e}rifi{\'e} pour toute tour de Bott de dimension $N-1$. Soit
$Y=Y_{C}$ et $Y'=Y_{C_{N-1}}$. Alors $Y=\mathbb{P}(\mathbf{1}\oplus \mathbf{L}_N)$, o{\`u} 
 $\mathbf{L}_N$ est un fibr{\'e} en droites au dessus de $Y'$ 
(voir la section~\ref{def} pour la d{\'e}finition de
$Y_{C_{N-1}}$ et  $\mathbf{L}_N$). On note 
$\pi_N : Y=\mathbb{P}(\mathbf{1}\oplus \mathbf{L}_N)
\rightarrow Y'$ la projection de $Y$ sur $Y'$ et $\mathbf{E}_N \in
K_D(Y)$ le fibr{\'e} tautologique au dessus de 
$Y=\mathbb{P}(\mathbf{1}\oplus \mathbf{L}_N)$.

\smallskip

Soit $D'=(S^1)^{N-1}$, par hypoth{\`e}se de r{\'e}currence, $K_{D'}(Y')$ est un 
$R[D']$-module libre engendr{\'e} par une base 
$\{ \hat{\mu}_{f}^{D'} \}_{f \in \{0,1\}^{N-1}}$. 
On d{\'e}finit une action de
 $D=(S^1)^N$
sur $Y'$ en faisant agir la derni{\`e}re composante trivialement sur $Y'$. On a alors 
$K_{D}(Y')=K_{D'}(Y')\otimes_{\mathbb{Z}}R[S^1]$ (voir \cite{kequivariante}),
et $K_{D}(Y')$ est donc
un $R[D]$-module libre qui admet comme base la famille
$\{ \hat{\mu}_{f}^D \}_{f \in \{0,1\}^{N-1}}$, o{\`u} pour tout $f \in
\{0,1\}^{N-1}$, on a pos{\'e} $ \hat{\mu}_{f}^D=\hat{\mu}_{f}^{D'}\otimes 1$.

\smallskip

 La projection $\pi_N$ est $D$-{\'e}quivariante et munit donc $K_D(Y)$
 d'une structure de $K_D(Y')$-module. Comme
$Y=\mathbb{P}(\mathbf{1}\oplus \mathbf{L}_N)$, 
 $K_{D}(Y)$ est un $K_{D}(Y')$-module 
libre de rang $2$ engendr{\'e} par le fibr{\'e} trivial $\mathbf{1}$ et le fibr{\'e} tautologique  
$\mathbf{E}_N$ soumis {\`a} l'unique relation  
\begin{equation} \label{EN2}
\mathbf{E}_N^2=(\mathbf{1}\oplus \mathbf{L}_N) \mathbf{E}_N -
\mathbf{L}_N \mathbf{1}, \end{equation}
(voir \cite{atiyah} th{\'e}or{\`e}me 2.2.1, ou \cite{kequivariante}). 
Soit $p : \mathcal{E}=\{0,1\}^N \rightarrow
\{0,1\}^{N-1}$ la projection selon les $N-1$ premi{\`e}res coordonn{\'e}es, 
on obtient donc une base du  
$R[D]$-module $K_D(Y)$ en posant

$$\hat{\mu}_{\epsilon}^D = \left\{ \begin{array}{ll}\vspace{0,1 cm}
\pi_N^*(\hat{\mu}_{p(\epsilon)}^D)\mathbf{E}_N
 &\,\,\, { \rm si } \,\,\, \epsilon_N =0, \\ 

\pi_N^*(\hat{\mu}_{p(\epsilon)}^D)(\mathbf{1}-\mathbf{E}_N)
 &\,\,\, { \rm si } \,\,\, \epsilon_N =1 .
\end{array}\right.$$
 
\end{proof}

On va expliciter la base $\{ \hat{\mu}_{\epsilon}^D \}_{\epsilon \in
  \mathcal{E}}$ d{\'e}finie dans la d{\'e}monstration de la proposition pr{\'e}c{\'e}dente. On
  reprend le diagramme de la section~\ref{def} : 
$$ \begin{array}{clccccc}
 & &  & & & \mathbb{P}(\mathbf{1}\oplus \mathbf{L}_N)=&Y_C \\ 
 &  &  & &  &  \downarrow \pi_N &  \\
  &  & & & &  Y_{C_{N-1}}&  \\
 & &   & &\Ddots &  \\
  & & \mathbb{P}(\mathbf{1}\oplus  \mathbf{L}_2)=& Y_{C_2}&   & &  \\
  & & \downarrow \pi_2 & &  & &  \\
  & \mathbb{C}P^1 \,\,\, = & \!\! Y_{C_1} & & & &  \\
  &\,\, \downarrow \pi_1  & & &  & &  \\
 \{un \, point\}= &\,\,  Y_0 & & &  & &  
\end{array}$$

\medskip

Pour $1 \leq i \leq N$, chaque vari{\'e}t{\'e} $Y_{C_i}$ (o{\`u} $Y_{C_N}=Y_C$)  
est munie de l'action de $D$
induite par l'action~\ref{actiontours} de $D_{\mathbb{C}}$, 
et les projections $\pi_i$ sont $D$-{\'e}quivariantes.
Pour $0 \leq i \leq N$, on pose $\Pi_i=\pi_{i+1} \pi_{i+2} \cdots \pi_N  : Y 
\rightarrow Y_{C_{i}}$ ($\Pi_N = { \rm Id }_{Y}$ et $Y_{C_{0}}=Y_0$). 

On pose $\mathbf{E}_i= \Pi_i^*(\mathbf{T}_i)\in K_{D}(Y)$, o{\`u} 
$\mathbf{T}_i \in K_{D}(Y_{C_i})$ est le fibr{\'e}
 tautologique sur $Y_{C_{i}}$, et 
 $\mathbf{F}_i=\mathbf{1}-\mathbf{E}_i $.

Pour tout $\epsilon \in \mathcal{E}$, la classe 
$\hat{\mu}_{\epsilon}^D \in K_D(Y)$ est alors donn{\'e}e par la formule
suivante :
\begin{equation} \label{muproduit}
 \boxed{ \hat{\mu}_{\epsilon}^D =\prod_{i\in \pi_+(\epsilon)} \!\!\!\! \mathbf{F}_i
\prod_{j\in \pi_-(\epsilon)}\!\!\!\!\mathbf{E}_j.}
\end{equation}

\subsection{Restrictions aux points fixes}

Pour tout $1 \leq i \leq N$, la projection $\Pi_{i-1} :  Y \rightarrow
Y_{C_{i-1}}$ permet de definir un {\'e}l{\'e}ment de $K_D(Y)$ {\`a} partir du
fibr{\'e} en droites  $D$-{\'e}quivariant $\mathbf{L}_i$ sur
$Y_{C_{i-1}}$. On note encore $\mathbf{L}_i$ cet {\'e}l{\'e}ment. 

Les d{\'e}finitions~\ref{actiontours} et \ref{actionfibres} de l'action de
$D$ sur $Y_{C_i}$ et sur $\mathbf{L}_i$ permettent de calculer
facilement les restrictions aux points fixes des fibr{\'e}s en droites 
$\mathbf{E}_i$ et $\mathbf{L}_i$ :
 
\begin{lemm} \label{toutesrestrictions}

$$i_D^*( \mathbf{E}_i)(\epsilon) = \left\{ \begin{array}{ll} 1

 & { \rm si } \hspace{0,15 cm} \epsilon_i = 0,  \\ 
e^{-\lambda_i(\epsilon)}
 & { \rm si } \hspace{0,15 cm} \epsilon_i = 1,
\end{array}\right.$$

$$i_D^*( \mathbf{F}_i)(\epsilon)=\left\{ \begin{array}{ll} 0

 & { \rm si } \hspace{0,15 cm} \epsilon_i = 0,  \\ 
1-e^{-\lambda_i(\epsilon)}
 & { \rm si } \hspace{0,15 cm} \epsilon_i = 1,
\end{array}\right.$$

$$ i_D^*( \mathbf{L}_i)(\epsilon)=
e^{-\lambda_i}\!\!\!\!\! \!\!\!\!\prod_{\tiny \begin{array}{cc}  1
    \leq j \leq i-1 \\ j \in 
    \pi_+(\epsilon) \end{array}} \!\!\!\!\! \!
e^{-c_{j,i}(\epsilon)\lambda_j} .$$

\end{lemm}





Le th{\'e}or{\`e}me suivant donne la valeur des restrictions aux points fixes des
classes $\hat{\mu}_{\epsilon}^D $.
Si on pose $ \mu_{\epsilon}^D = i_{D}^*( \hat{\mu}_{\epsilon}^D)$, le
lemme~\ref{toutesrestrictions} et la formule~\ref{muproduit} donnent
imm{\'e}diatement la formule suivante :  
\begin{theo}  \label{restrictionkb}
    
Pour $(\epsilon,\epsilon') \in \mathcal{E}^2$,

$$ \mu_{\epsilon}^D(\epsilon') = \left\{ \begin{array}{ll}
\!\!\!\! \displaystyle{\prod_{i \in \pi_{+}(\epsilon')}e^{-\lambda_{i}(\epsilon')}
\prod_{i \in \pi_{+}(\epsilon)}(e^{\lambda_{i}(\epsilon')}-1) }
 & { \rm si } \hspace{0,15 cm} \epsilon \leq \epsilon', \\ 
  0
 & { \rm  sinon }.
\end{array}\right.$$

\end{theo}

\begin{exem}

On consid{\`e}re la surface de Hirzebruch $H_{-1}=Y_{\{-1\}}$. On pose 

 $$\left\{ \begin{array}{lll} \epsilon^{1}=(0,0), &   \\ 
 \epsilon^{2}=(1,0), & \epsilon^{3}=(0,1) ,
 \\ \epsilon^{4}=(1,1). & 
 & \end{array} \right.$$

Si on d{\'e}finit la matrice $\mathcal{M}=\{ \mu_{i,j}  \}_{1 \leq i< j \leq 4}$ par 
$\mu_{i,j}=\mu_{\epsilon^i}^D(\epsilon^j)$, on obtient les formules
suivantes :

$$\mathcal{M}=\begin{pmatrix} 1 & e^{-\lambda_1} & e^{-\lambda_2} &
  e^{-2\lambda_1-\lambda_2} \\
 0 &1- e^{-\lambda_1} & 0 &
  e^{-\lambda_1-\lambda_2}(1- e^{-\lambda_1}) \\
 0 & 0 &1- e^{-\lambda_2} &
  e^{-\lambda_1}(1- e^{-\lambda_1-\lambda_2}) \\
0 & 0 &0  &
  (1- e^{-\lambda_1})(1-e^{-\lambda_1-\lambda_2}) 
\end{pmatrix}.$$

\end{exem}

\subsection{Caract{\'e}risation de la base $\hat{\mu}_{\epsilon}^D$}

La base $\{  \hat{\mu}_{\epsilon}^D \} _{\epsilon \in \mathcal{E}}$ est reli{\'e}e {\`a}
la d{\'e}composition cellulaire $Y=\coprod_{\epsilon \in
  \mathcal{E}}Y_{\epsilon}$ par le th{\'e}or{\`e}me suivant :

\begin{theo} \label{ktheorietours}
    
La famille $\{  \hat{\mu}_{\epsilon}^D \} _{\epsilon \in \mathcal{E}}$ est une
base du R[D]-module $K_{D}(Y)$ caract{\'e}ris{\'e}e par les relations 
$$\forall
(\epsilon, \epsilon') \in \mathcal{E}^2, \,\, \chi(\overline{Y}_{\epsilon'}, 
\hat{\mu}_{\epsilon}^D)= \delta_{\epsilon, \epsilon'} . $$
   
\end{theo}

\begin{proof}

On sait d{\'e}j{\`a} que la famille  $\{  \hat{\mu}_{\epsilon}^D \} _{\epsilon \in
  \mathcal{E}}$ est une base de $K_{D}(Y)$. Pour $(\epsilon, \epsilon') \in
  \mathcal{E}^2$, on va calculer  
$\chi(\overline{Y}_{\epsilon'}, \hat{\mu}_{\epsilon}^D)$ gr{\^a}ce {\`a} la formule de
  localisation.

En utilisant la proposition~\ref{pointsfixesab}, 
 et le lemme~\ref{pointsfixestours}, on
obtient pour tout  $\hat{\mu}^D \in K_{D}(Y)$ et tout $\epsilon \in \mathcal{E}$: 
\begin{equation} \label{ab} \chi(\overline{Y}_{\epsilon}, \hat{\mu}^D)=  
\sum_{\epsilon' \leq
\epsilon}\frac{i_{D}^*(\hat{\mu}^D)(\epsilon')}{\prod_{i \in \pi_{+}(\epsilon)}
(1-e^{-\lambda_{i}(\epsilon')})}.  \end{equation}

Cette formule et le th{\'e}or{\`e}me~\ref{restrictionkb} nous montrent imm{\'e}diatement
que $\chi(\overline{Y}_{\epsilon}, \hat{\mu}_{\epsilon}^D)=1$ et 
$\chi(\overline{Y}_{\epsilon}, \hat{\mu}_{\epsilon_0}^D)=0$ si $\epsilon_0 \not\leq
\epsilon$.

Soit $\epsilon_0 \in \mathcal{E}$ tel que $\epsilon_0 \leq \epsilon$ et $\epsilon_0 \neq
\epsilon$. Alors, la formule~\ref{ab} et le th{\'e}or{\`e}me~\ref{restrictionkb} nous
donnent 
$$ \chi(\overline{Y}_{\epsilon}, \hat{\mu}_{\epsilon_0}^D) = \sum_{\epsilon_{0} \leq
  \epsilon' \leq
\epsilon}\frac{\displaystyle{\prod_{i \in \pi_{+}(\epsilon')}e^{-\lambda_{i}(\epsilon')}
\prod_{i \in \pi_{+}(\epsilon_{0})}(e^{\lambda_{i}(\epsilon')}-1)}}
{\displaystyle{\prod_{i \in
  \pi_{+}(\epsilon)}(1-e^{-\lambda_{i}(\epsilon')})}},$$
d'o{\`u} l'on tire 
$$\chi(\overline{Y}_{\epsilon}, \hat{\mu}_{\epsilon_0}^D) = 
\sum_{\epsilon_{0} \leq  \epsilon' \leq
\epsilon}\frac{\displaystyle{\prod_{i \in \pi_{+}(\epsilon') \setminus \pi_{+}(\epsilon_{0})}
e^{-\lambda_{i}(\epsilon')}} }
{\displaystyle{\prod_{i \in \pi_{+}(\epsilon)\setminus \pi_{+}(\epsilon_{0})}
(1-e^{-\lambda_{i}(\epsilon')})}}.$$

\medskip

Soit $j$ le plus grand
{\'e}l{\'e}ment de $\pi_{+}(\epsilon) \setminus  \pi_{+}(\epsilon_{0})$, on a alors 
$$\mbox{ \small { \mbox{ $ \chi(\overline{Y}_{\epsilon},
      \hat{\mu}_{\epsilon_0}^D)=$}}}\!\!\!\!\!\sum_{
 \tiny \begin{array}{cc}  \epsilon_{0} \leq \epsilon' \leq
\epsilon \\
 \epsilon'_j=0 
\end{array}}\!\!\!\!\!\!\!\!\mbox{ \small { \mbox{ $ 
\frac{ \displaystyle{ \prod_{i \in \pi_{+}(\epsilon') \setminus
  \pi_{+}(\epsilon_{0})} e^{-\lambda_{i}(\epsilon') }}}
{\displaystyle{\prod_{i \in \pi_{+}(\epsilon)\setminus \pi_{+}(\epsilon_{0})}
(1-e^{-\lambda_{i}(\epsilon')})}}$}}} +
\!\!\!\!\!\sum_{
 \tiny \begin{array}{cc}  \epsilon_{0} \leq \epsilon' \leq
\epsilon \\
 \epsilon'_j=1 
\end{array}}\!\!\!\!\!\!\!\!\mbox{ \small { \mbox{ $ 
\frac{ \displaystyle{ \prod_{i \in \pi_{+}(\epsilon') \setminus
  \pi_{+}(\epsilon_{0})} e^{-\lambda_{i}(\epsilon') }}}{\displaystyle{\prod_{i 
\in \pi_{+}(\epsilon)\setminus \pi_{+}(\epsilon_{0})}
(1-e^{-\lambda_{i}(\epsilon')})}}$}}}, $$
d'o{\`u} l'on tire 
$$\mbox{ \small { \mbox{ $ \chi(\overline{Y}_{\epsilon},
      \hat{\mu}_{\epsilon_0}^D)=$}}}\!\!\!\!\!\sum_{
 \tiny \begin{array}{cc}  \epsilon_{0} \leq \epsilon' \leq
\epsilon \\
 \epsilon'_j=0 
\end{array}}\!\!\!\!\!\!\!\! \mbox{ \small { \mbox{ $ \left[ \frac{\displaystyle{ 
\prod_{i \in \pi_{+}(\epsilon') \setminus
  \pi_{+}(\epsilon_{0})} e^{-\lambda_{i}(\epsilon') }}}{
\displaystyle{\prod_{i \in \pi_{+}(\epsilon)\setminus \pi_{+}(\epsilon_{0})}
(1-e^{-\lambda_{i}(\epsilon')})}}+\frac{ \displaystyle{ \prod_{i \in
  \pi_{+}(\epsilon'+(j)) \setminus
  \pi_{+}(\epsilon_{0})} e^{-\lambda_{i}(\epsilon'+(j)) }}}
{\displaystyle{\prod_{i \in \pi_{+}(\epsilon)\setminus \pi_{+}(\epsilon_{0})}
(1-e^{-\lambda_{i}(\epsilon'+(j))})}} \right]$}}}.
$$

Chaque terme de cette somme est nulle.
En effet, soit $\epsilon'$ un {\'e}l{\'e}ment de la sommation, 
comme $j$ est le plus grand {\'e}l{\'e}ment de $\pi_{+}(\epsilon) \setminus
\pi_{+}(\epsilon_{0})$, pour tout $i \in \pi_{+}(\epsilon) \setminus
\pi_{+}(\epsilon_{0})$, 
$\lambda_{i}(\epsilon'+(j))=\lambda_{i}(\epsilon')$ si $i\neq j$, et 
$\lambda_{j}(\epsilon'+(j))=-\lambda_{j}(\epsilon')$. Le terme de la somme
associ{\'e} {\`a} $\epsilon'$ est donc 
$$\frac{ \displaystyle{ \prod_{i \in \pi_{+}(\epsilon') \setminus
  \pi_{+}(\epsilon_{0})} e^{-\lambda_{i}(\epsilon') }}}{\displaystyle{\prod_{i 
\in \pi_{+}(\epsilon-(j))\setminus \pi_{+}(\epsilon_{0})}
(1-e^{-\lambda_{i}(\epsilon')})}} \,\,\,\, \left[ \frac{1} 
{1-e^{-\lambda_{j}(\epsilon')} }+ 
 \frac{e^{\lambda_{j}(\epsilon')}  }{ 1-e^{\lambda_{j}(\epsilon')} }
 \right].$$
 
Ce terme est bien nul d'apr{\`e}s la relation 
$\frac{1}{1-e^{-x}}+\frac{e^{x}}{1-e^{x}}=0$, et on obtient donc $\chi(
\overline{Y}_{\epsilon}, \hat{\mu}_{\epsilon_0}^D) =0$.

\end{proof}

\subsection{Structure multiplicative} \label{multiplicationtours}

La $R[D]$-alg{\`e}bre $K_D(Y)$ est engendr{\'e}e par les $\mathbf{E}_i$. On
rappelle qu'on a pos{\'e} $\mathbf{F}_i=\mathbf{1}-  \mathbf{E}_i$ et que
la famille $\{  \hat{\mu}_{\epsilon}^D \}_{\epsilon \in
  \mathcal{E}}$ d{\'e}finie par $$\hat{\mu}_{\epsilon}^D = \prod_{i\in
  \pi_+(\epsilon)} \!\!\!\! \mathbf{F}_i
\prod_{j\in \pi_-(\epsilon)}\!\!\!\!\mathbf{E}_j$$ est une base du
$R[D]$-module  $K_{D}(Y)$.

\subsubsection{Pr{\'e}sentation par g{\'e}n{\'e}rateurs et relations}

Pour tout $1 \leq i \leq N$, la relation~\ref{EN2} s'{\'e}crit sous la forme : 
\begin{equation}\label{Ei2}
\mathbf{E}_i^2=(\mathbf{1} + \mathbf{L}_i) \mathbf{E}_i -
\mathbf{L}_i =\mathbf{E}_i  -\mathbf{L}_i\mathbf{F}_i .
\end{equation}

On en d{\'e}duit facilement les relations suivantes : 
\begin{equation} \label{relations}
\left\{
\begin{array}{c}
\mathbf{F}_i^2=(\mathbf{1}-\mathbf{L}_i) \mathbf{F}_i ,\\
\mathbf{E}_i\mathbf{F}_i=\mathbf{L}_i\mathbf{F}_i,\\
\mathbf{E}_i^{-1}=*\mathbf{E}_i=\mathbf{E}_i+(\mathbf{1}+\mathbf{L}_i^{-1})\mathbf{F}_i.
\end{array}  \right.
\end{equation}

Pour comprendre la structure multiplicative de $K_D(Y)$, il faut donc
exprimer $\mathbf{L}_i$ en fonction des  $\mathbf{E}_j$ pour $1  \leq
j < i$ :

\begin{lemm} \label{li}

\indent

$$\mathbf{L}_i = e^{-\lambda_i}\prod_{1\leq
  j<i}\mathbf{E}_{j}^{-c_{j,i}}.
$$

\end{lemm}

\begin{proof}

On peut d{\'e}montrer ce lemme gr{\^a}ce aux restrictions aux points
fixes. Comme l'application $i_D^* : K_D(Y) \rightarrow F(\mathcal{E};
R[D])$ est un morphisme  de $R[D]$-alg{\`e}bres
injectif, il suffit de montrer
que pour tout $\epsilon \in \mathcal{E}$,
$$i_D^*(\mathbf{L}_i)(\epsilon)
=e^{-\lambda_i}\prod_{1\leq
  j<i}\Big(i_D^*(\mathbf{E}_{j})(\epsilon)\Big)^{-c_{j,i}}.$$

D'apr{\`e}s le lemme~\ref{toutesrestrictions}, on doit donc montrer
que pour tout $\epsilon \in \mathcal{E}$,
$$\sum_{\tiny \begin{array}{c} 1 \leq j <i \\ j \in
    \pi_+(\epsilon) \end{array}}\!\!\!\! -c_{j,i}(\epsilon)\lambda_j=
\!\!\!\! \sum_{\tiny \begin{array}{c} 1 \leq j <i \\ j \in
    \pi_+(\epsilon) \end{array}}\!\!\!\! c_{j,i}\lambda_j(\epsilon).$$
 
Dans la deuxi{\`e}me somme ci-dessus, le coefficient de $\lambda_j$ est
{\'e}gal {\`a}  $$c_{j,i} + \!\!\!\! \!\!\sum_{\tiny \begin{array}{c} j< k < i  \\ k \in
    \pi_+(\epsilon) \end{array}} \!\!\!\!
\!\!c_{k,i}c_{j,k}(\epsilon)= -c_{j,i}(\epsilon),$$
d'apr{\`e}s la relation~\ref{c}, et on a donc bien l'{\'e}galit{\'e} voulue.

\end{proof} 

On note $\mathcal{B}$ la $R[D]$-alg{\`e}bre $R[D][X_1, \ldots , X_N, X_1^{-1},
 \ldots ,  X_N^{-1}]$, et pour $1 \leq i \leq N$, on d{\'e}finit
l'{\'e}l{\'e}ment $L_i$ de $\mathcal{B}$ par 
  $$\boxed{L_i = e^{-\lambda_i}\prod_{1 \leq j
        <i}X_j^{-c_{j,i}}.}$$

On peut alors r{\'e}sumer les r{\'e}sultats pr{\'e}c{\'e}dents dans le th{\'e}or{\`e}me
suivant :

\begin{theo} \label{generateursrelationstours}
 Le morphisme de $R[D]$-alg{\`e}bres de $\mathcal{B}$
dans $K_D(Y)$ qui envoie $X_i$ sur $\mathbf{E}_i$ et $X_i^{-1}$ sur
$\mathbf{E}_i^{-1}=*\mathbf{E}_i$  induit un
 isomorphisme $$ f :  \mathcal{C}=\mathcal{B}/  \mathcal{I} \rightarrow  K_D(Y),$$ o{\`u}
 $\mathcal{I}$ est  l'id{\'e}al engendr{\'e} par
$\{ X_i^2-X_i+(1-X_i)L_i \}_{1
  \leq i \leq N}$.

 \smallskip

 De plus, la famille de polyn{\^o}mes $\{ Q_{\epsilon} \}_{\epsilon \in
   \mathcal{E}}$ d{\'e}finie par  $$ Q_{\epsilon}=\prod_{i \in
    \pi_-(\epsilon)}\!\!\! X_i \prod_{j \in
   \pi_+(\epsilon)}\!\!(1-X_j) \in \mathcal{C}{ \,\,\,\,\,\,\,\, \rm  \,
   pour \,  tout \, }  \epsilon \in
   \mathcal{E},$$ est une base du $R[D]$-module libre $\mathcal{C}$, et
   $Q_{\epsilon}$ est envoy{\'e} sur $\hat{\mu}_{\epsilon}^D$ par
   l'isomorphisme~$f$.

 \end{theo}

 \subsubsection{Calcul des constantes de structure}

 On va donner une m{\'e}thode de calcul des {\'e}l{\'e}ments $r_{\epsilon, \epsilon
  '}^{\epsilon''} \in R[D]$ qui v{\'e}rifient les relations 
$$\hat{\mu}_{\epsilon}^D\hat{\mu}_{\epsilon'}^D=
\sum_{\epsilon'' \in \mathcal{E}} r_{\epsilon, \epsilon
  '}^{\epsilon''}\hat{\mu}_{\epsilon''}^D.
$$

Ces {\'e}l{\'e}ments sont {\'e}galement d{\'e}finis par les relations
$$  r_{\epsilon, \epsilon'}^{\epsilon''}=
 \chi(\overline{Y}_{\epsilon''},
 \hat{\mu}_{\epsilon}^D\hat{\mu}_{\epsilon'}^D).
$$

\bigskip

On introduit la $R[D]$-alg{\`e}bre
$$\mathcal{D}=\mathcal{B}[Z_1, \ldots Z_N]=R[D][ X_1, \ldots ,  X_N, 
X_1^{-1}, \ldots  , X_N^{-1}, Z_1, \ldots , Z_N].$$

On note $\beta  :
\mathcal{D}  \rightarrow \mathcal{C}$ le morphisme
induit par le morphisme de $\mathcal{B}$-alg{\`e}bres de $\mathcal{D}$
dans $\mathcal{B}$ qui envoie $Z_i$ sur $1-X_i$.

\bigskip

\begin{defi}
Pour tout $1 \leq i \leq N+1$, on note $o_i : \mathcal{D} \rightarrow
\mathcal{D} $ le morphisme de $R[D]$-alg{\`e}bres d{\'e}fini par 
$$ o_i(X_j)= \left\{ \begin{array}{ll}1  & {\rm \, 
pour \, } j \geq i, \\
X_j & {\rm \, 
pour \, } j < i,
\end{array} \right. \,\,\,\,\,\,\,\,\,
 o_i(Z_j)= \left\{ \begin{array}{ll} 0 & {\rm \, 
pour \, } j \geq i, \\
Z_j & {\rm \, 
pour \, } j < i.
\end{array} \right.
$$

L'application $o_{N+1}$ est l'identit{\'e} de $\mathcal{D}$.

\end{defi}

\begin{defi}
On note $\omega_i : \mathcal{B} \rightarrow \mathcal{B}$ le morphisme de 
$R[D]$-alg{\`e}bres d{\'e}fini par 
 $$\omega_i(X_j)= \left\{ \begin{array}{ll} 1 & {\rm \, 
pour \, } j \geq i, \\
X_j &  {\rm \, 
pour \, } j < i,
\end{array} \right.
$$
et on note $\overline{\omega}_i$ le morphisme induit sur $\mathcal{C}$
(ce morphisme est bien d{\'e}fini car $\omega_i$ laisse stable
$\mathcal{I}$). 
\end{defi}

\begin{exem}\label{omegaQ}
$$ \overline{w}_i(Q_{\epsilon} )= \left\{ \begin{array}{ll} 
 Q_{\epsilon}X_i^{-1} \cdots X_N^{-1}  &  {\rm \, 
si  \,\, } \forall j \geq i , \, \epsilon_j=0 \\
0 & {\rm \, 
sinon \, }.
\end{array} \right.
$$

\end{exem} 

\smallskip

On a la relation  \begin{equation} \label{commutation}
\beta \circ o_i = \overline{\omega}_i \circ \beta.
\end{equation}

\begin{defi} \label{defR}

Soit $\epsilon$ un {\'e}l{\'e}ment de  $\mathcal{E}$ de longueur $l$ strictement positive. 
On note $\{i_1 < \cdots <i_l \}$ les {\'e}l{\'e}ments
 de $\pi_+(\epsilon)$. On d{\'e}finit alors l'application $R^{\epsilon} : 
\mathcal{D} \rightarrow R[D]$ de la mani{\`e}re suivante :

\begin{enumerate}

\item[$(i)$] $R^{\epsilon}$ est $R[D]$-lin{\'e}aire,


\item[$(ii)$] si $P\in \mathcal{D}$ est un mon{\^o}me non nul qui s'{\'e}crit 
sous la forme $P=SX_{i_l}^r Z_{i_l}^s$
  o{\`u} $S \in \mathcal{D}$ est de degr{\'e} $0$ en  $X_{i_l}$
   et $Z_{i_l}$, et o{\`u} $s$ est un entier positif et $r$
  un entier quelconque, 
$$R^{\epsilon}(P)= \left\{ \begin{array}{ll}
R^{\epsilon-(i_l)}\Big(S(1-L_{i_l})^{s-1}L_{i_l}^r\Big) & { \rm si \,
} s>0, \\
-R^{\epsilon-(i_l)}\Big(S(L_{i_l}+L_{i_l}^2 + \cdots + L_{i_l}^{r-1})\Big)
& { \rm si \,} s=0 { \rm  \, \, et \, } r>1 ,\\
R^{\epsilon-(i_l)}\Big(S(1+L_{i_l}^{-1} + \cdots
+ L_{i_l}^{r})\Big)
& { \rm si \,} s=0 { \rm \, \,  et \, } r < 0, \\
R^{\epsilon-(i_l)}(S) & { \rm si \,} s=0 { \rm \, \,
  et \, } r = 0 ,\\
0 & { \rm si \,} s=0 { \rm \, \,
  et \, } r = 1,
\end{array} \right.
$$
\item[$(iii)$] $R^{(\mathbf{0})}(P)=P(X_i=1, Z_i=0)$.

\end{enumerate}

\smallskip

Ces trois relations d{\'e}finissent compl{\`e}tement (r{\'e}cursivement) les applications
 $R^{\epsilon}$.

\end{defi}

\begin{rema}

La d{\'e}finition de ces applications $R^{\epsilon}$ est inspir{\'e}e de celle des
applications $T_A$ d{\'e}finies par Haibao Duan en cohomologie ordinaire
dans \cite{duan3} et g{\'e}n{\'e}ralis{\'e}es en cohomologie {\'e}quivariante
dans \cite{mw4}.

\end{rema}

\begin{theo} \label{theoP}

 Pour $\epsilon \in
\mathcal{E}$, on note $\beta_{\epsilon} : \mathcal{D} \rightarrow
R[D]$ l'application d{\'e}termin{\'e}e
par les relations  \begin{equation} \label{betaepsilon}
\forall P \in \mathcal{D}, \,\, \beta(P)=\sum_{\epsilon \in
   \mathcal{E}}\beta_{\epsilon}(P)Q_{\epsilon}. 
\end{equation}

 Alors, pour tout $\epsilon \in \mathcal{E}$, 
$\beta_{\epsilon}=R^{\epsilon}$.
 
 \end{theo}

On en d{\'e}duit imm{\'e}diatement le corollaire suivant : 

\begin{coro} \label{coror}

Pour $\epsilon \in \mathcal{E}$, on pose $S_{\epsilon}^D=\prod_{i \in
   \pi_-(\epsilon)} X_i \prod_{j \in \pi_+(\epsilon)}Z_j \in
   \mathcal{D}$.
Soient $\epsilon$, $\epsilon'$ et $\epsilon''$ trois {\'e}l{\'e}ments de
$\mathcal{E}$, alors : 
 $$r_{\epsilon, \epsilon'}^{\epsilon''}=R^{\epsilon''}
(S_{\epsilon}^D S_{\epsilon'}^D).$$

\end{coro}

\begin{rema}

On montre facilement {\`a} l'aide de ce corollaire que $r_{\epsilon,
  \epsilon'}^{\epsilon''}=0$ sauf si $\epsilon'' \geq \epsilon$ et 
$\epsilon'' \geq \epsilon'$.

\end{rema}

Pour d{\'e}montrer le th{\'e}or{\`e}me~\ref{theoP}, on aura besoin du lemme suivant : 

\begin{lemm} \label{beta}

Soit $\epsilon$ un {\'e}l{\'e}ment de  $\mathcal{E}$ de longueur $l$ strictement positive. 
On note $\{i_1 < \cdots <i_l \}$ les {\'e}l{\'e}ments
 de $\pi_+(\epsilon)$. Si $S$ est un mon{\^o}me non nul de $\mathcal{D}$ de degr{\'e}
 nul en $X_{i_l}$ et $Z_{i_l}$, alors :  
$$\left\{ \begin{array}{lll} \beta_{\epsilon}(SX_{i_l}) & = &0, \\
\beta_{\epsilon}(SZ_{i_l}) & =  & \beta_{\epsilon-(i_l)}(S), \\
\beta_{\epsilon}(S) & =  & \beta_{\epsilon-(i_l)}(S).
\end{array} \right.
$$

\end{lemm}

\begin{proof}

Comme $S$ est un mon{\^o}me non nul de degr{\'e} $0$ en $X_{i_l}$ et
$Z_{i_l}$, $$o_{i_l}(S)=o_{i_l+1}(S).$$

On a donc  \begin{equation} \label{formule}
\beta \big(o_{i_l+1}(SZ_{i_l}) \big)=\beta \big(o_{i_l+1}(S)
\big)\beta \big(o_{i_l+1}(Z_{i_l})\big)= \beta
\big(o_{i_l}(S)\big)(1-X_{i_l}).
\end{equation}

Le premier membre de l'{\'e}galit{\'e}~\ref{formule} est {\'e}gal {\`a} 
$$\beta \big(o_{i+1}(SZ_{i_l})
\big)=\overline{\omega}_{i_l+1}\big( \beta(SZ_{i_l}) \big)=\!\!\!\!\!\!\sum_{\tiny
  \begin{array}{c} \epsilon' \in \mathcal{E} \\ \epsilon'_j=0 \, \, \forall j
    > i_l \end{array} } \!\!\!\!\!\!\beta_{\epsilon'} (SZ_{i_l}) 
Q_{\epsilon'}X_{i_l+1}^{-1} \cdots X_{N}^{-1}.$$

   La premi{\`e}re {\'e}galit{\'e} provient de la relation de
  commutation~\ref{commutation} et la seconde de l'exemple~\ref{omegaQ}.

  \smallskip

  Le dernier membre de l'{\'e}galit{\'e}~\ref{formule} est {\'e}gal {\`a} 
$$\beta \big(o_i(S)\big)(1-X_{i_l})
=\overline{\omega}_{i_l}\big( \beta(S) \big)(1-X_{i_l})=\!\!\!\!\!\!\sum_{\tiny
  \begin{array}{c} \epsilon'' \in \mathcal{E} \\ \epsilon''_j=0 \, \, \forall j
    \geq i_l \end{array} } \!\!\!\!\!\!\beta_{\epsilon''} (S) 
Q_{\epsilon''}(1-X_{i_l}) X_{i_l+1}^{-1} \cdots X_{N}^{-1}$$
$$=\sum_{\tiny
  \begin{array}{c} \epsilon' \in \mathcal{E}, \,  \epsilon'_{i_l}=1 
\\ \epsilon'_j=0 \, \, \forall j
    > i_l \end{array} } \!\!\!\!\!\!\beta_{\epsilon'-(i_l)} (S) 
Q_{\epsilon'} X_{i_l+1}^{-1} \cdots X_{N}^{-1}.
$$

On en d{\'e}duit donc l'{\'e}galit{\'e} : 
$$  \sum_{\tiny
  \begin{array}{c} \epsilon' \in \mathcal{E} \\ \epsilon'_j=0 \, \, \forall j
    > i_l \end{array} } \!\!\!\!\!\!\beta_{\epsilon'} (SZ_{i_l}) 
Q_{\epsilon'}=\!\!\!\!\!\!
\sum_{\tiny
  \begin{array}{c} \epsilon' \in \mathcal{E}, \,  \epsilon'_{i_l}=1 
\\ \epsilon'_j=0 \, \, \forall j
    > i_l \end{array} } \!\!\!\!\!\!\beta_{\epsilon'-(i_l)} (S) 
Q_{\epsilon'}.
$$

La famille $\{ Q_{\epsilon} \}_{\epsilon \in \mathcal{E}}$ {\'e}tant une
base de $\mathcal{C}$ sur $R[D]$, on en d{\'e}duit la deuxi{\`e}me {\'e}galit{\'e} du
lemme. La premi{\`e}re se d{\'e}montre exactement de la m{\^e}me mani{\`e}re, et la
troisi{\`e}me est une cons{\'e}quence des deux premi{\`e}res car
$\beta(S(X_{i_l}+Z_{i_l}))=\beta(S)$.

\end{proof}

\begin{proof}[D{\'e}monstration du th{\'e}or{\`e}me~\ref{theoP}]

Pour d{\'e}montrer le th{\'e}or{\`e}me, il faut montrer que les applications
$\beta_{\epsilon}$ v{\'e}rifient les trois relations de la
d{\'e}finition~\ref{defR}.

La relation $(i)$ est imm{\'e}diate.

Pour trouver la relation $(iii)$, il suffit d'appliquer
$\overline{\omega}_1$ aux 
formules~\ref{betaepsilon}.

Pour montrer la relation $(ii)$, on consid{\`e}re un {\'e}l{\'e}ment de
$\mathcal{E}$ de longueur $l$ strictement positive, et on note $\{i_1
< \cdots <i_l \}$  les {\'e}l{\'e}ments
 de $\pi_+(\epsilon)$.

Les relations~\ref{Ei2} et \ref{relations} permettent de montrer
facilement par r{\'e}currence les relations suivantes valables dans 
$\mathcal{C}$ pour tout $1 \leq i \leq N$ et tout entier relatif $n$, 
$$\left\{
\begin{array}{ll}
X_i^n = X_i - (L_i + L_i^2 + \cdots + L_i^{n-1})(1-X_i) & { \rm si \,
} n > 0, \\
X_i^n = X_i + (1+L_i^{-1}  + \cdots + L_i^{n})(1-X_i) & { \rm
  si \, } n  
< 0, \\
(1-X_i)^n = (1-L_i )^{n-1}(1-X_i)   & { \rm si \, } n >0, \\
X_i^n(1-X_i)=L_i^n(1-X_i)  &  \forall \,   n .
\end{array}\right.
$$

Supposons par exemple que $P \in \mathcal{D}$ est un mon{\^o}me non nul 
qui s'{\'e}crit sous la forme $P=SX_{i_l}^r $
  o{\`u} $S \in \mathcal{D}$ est de degr{\'e} $0$ en  $X_{i_l}$ et $Z_{i_l}$,
  et o{\`u} $r < 0$. Alors, d'apr{\`e}s les relations ci-dessus,
$$ \beta\big(SX_{i_l}^r \big)=\beta \Big(  S 
\big( X_{i_l} + (1+L_{i_l}^{-1}  + \cdots +
L_{i_l}^{r})Z_{i_l} \big)
\Big),
$$
et on a donc en particulier 
$$\beta_{\epsilon}(P)=\beta_{\epsilon}(SX_{i_l})+\beta_{\epsilon}\big(S
(1+L_{i_l}^{-1}  + \cdots +
L_{i_l}^{r})Z_{i_l}\big).
$$

Comme chaque terme $S$ et $SL_{i_l}^{-k}$ pour $1\leq k \leq -r$
 est un mon{\^o}me non nul
de degr{\'e} $0$ en $X_{i_l}$ et $Z_{i_l}$, on peut appliquer le
lemme~\ref{beta} pour obtenir  
$$
\beta_{\epsilon}(P)=\beta_{\epsilon-(i_l)}\big( S
(1+L_{i_l}^{-1} + \cdots +
L_{i_l}^{r})\big).
$$

Les autres cas se traitent de la m{\^e}me mani{\`e}re.

\end{proof}

\section{K-th{\'e}orie {\'e}quivariante des vari{\'e}t{\'e}s de Bott-Samelson} \label{sectionBott}

 On reprend les notations de la section~\ref{BS}. On choisit $N$ racines simples
$\mu_1, \ldots , \mu_N$ non n{\'e}cessairement distinctes, et on pose $\Gamma = 
\Gamma( \mu_1, \ldots , \mu_N  )$. Pour $1 \leq i<j \leq N$, on pose
$b_{i,j}= \mu_j(\mu^{\vee}_i)$ et $B=\{b_{i,j}\}_{1\leq i<j \leq N}$.

 Les r{\'e}sultats de cette section ne d{\'e}pendent
 pas de la structure complexe de $\Gamma$, et on identifie donc
$\Gamma$
 avec la tour de Bott $Y=Y_B$, la d{\'e}composition 
$\Gamma=\coprod_{\epsilon \in \mathcal{E}} \Gamma_{\epsilon}$ avec la
  d{\'e}composition $Y=\coprod_{\epsilon \in \mathcal{E}} Y_{\epsilon}$, et le point
  fixe $\epsilon \in \Gamma^T$ avec le point fixe $\epsilon \in Y^D$.

Le tore $T$  agit sur $\Gamma$ via son image  $S$, et l'action 
 de $S$ sur $\Gamma$ s'identifie {\`a} celle d'un sous-tore de
$D$ sur $Y$. Pour tout $\epsilon \in \mathcal{E}$, on note alors
$\hat{\mu}^S_{\epsilon}=p_S^D(\hat{\mu}^D_{\epsilon})$ l'{\'e}l{\'e}ment de
 $K_S(\Gamma)$ obtenu {\`a} partir de 
$\hat{\mu}^D_{\epsilon}$ par restriction {\`a} $S$ de l'action de $D$. Ces {\'e}l{\'e}ments
s'obtiennent {\`a} l'aide des fibr{\'e}s de Hopf de la m{\^e}me mani{\`e}re que les classes 
$\hat{\mu}^D_{\epsilon}$, et ils forment donc une base du $R[S]$-module 
$K_S(\Gamma)$.

De plus, pour tout couple $(\epsilon,\epsilon') \in \mathcal{E}^2$,
$$\chi(\overline{\Gamma}_{\epsilon'}, \hat{\mu}_{\epsilon}^S)
=\delta_{\epsilon', \epsilon}.
$$

Comme le tore $T$ agit sur $\Gamma$ via son image $S$, on a $R[S] \subset
R[T]$, et  
$K_T(\Gamma)$ s'identifie {\`a} 
$K_S(\Gamma)\! \otimes_{R[S]}\!\!R[T]$. Si on pose 
$\hat{\mu}^T_{\epsilon}=\hat{\mu}^S_{\epsilon}\otimes 1 $, on a donc la proposition
suivante :

\begin{prop} \label{propbasekbs}

\indent

\begin{enumerate}

\item[$(i)$] La  $K$-th{\'e}orie $T$-{\'e}quivariante de $\Gamma^T$ s'identifie {\`a}
$F(\mathcal{E};R[T])$.

\item[$(ii)$] La restriction aux points fixes $i_{T}^*$ : $K_{T}(\Gamma)
 \rightarrow F(\mathcal{E};R[T])$ est injective.
    
\item[$(iii)$] La $K$-th{\'e}orie $T$-{\'e}quivariante de $\Gamma$ est un $R[T]$-module libre
 admettant comme base la famille $\{\hat{\mu}_{\epsilon}^T\}_{\epsilon \in
\mathcal{E}}$ qui v{\'e}rifie les relations
$$\forall (\epsilon, \epsilon') \in \mathcal{E}^2, \,\,
\chi(\overline{\Gamma}_{\epsilon'}, \hat{\mu}_{\epsilon}^T)=\delta_{\epsilon',
  \epsilon}.$$
   
\end{enumerate}

\end{prop}

\subsection{Restrictions aux points fixes}

Pour tout $\epsilon \in \mathcal{E}$, on pose  $\mu_{\epsilon}^T=
i_{T}^*(\hat{\mu}_{\epsilon}^T)$. Le th{\'e}or{\`e}me~\ref{restrictionkb} nous
permet alors de retrouver la valeur des restrictions aux points fixes
de ces classes $\hat{\mu}_{\epsilon}^T$ calcul{\'e}es dans \cite{mw3} :

\begin{theo} \label{restrictionskbs}
    
Pour $\epsilon \in \mathcal{E}$,

$$\mu_{\epsilon}^T(\epsilon') = \left\{ \begin{array}{ll}  
\!\!\!\!\displaystyle{\prod_{i \in \pi_{+}(\epsilon')}e^{\alpha_{i}(\epsilon')}
\prod_{i \in \pi_{+}(\epsilon)}(e^{-\alpha_{i}(\epsilon')}-1) }
 & { \rm si } \hspace{0,15 cm} \epsilon \leq \epsilon', \\ 
 0
 & { \rm  sinon }.
\end{array}\right.$$
   
 \end{theo}

  \begin{proof}

 On a le diagramme commutatif suivant :

 $$\xymatrix{
     K_{S}(\Gamma)\ar@{^{(}->}[dd]^{i_{S}^*} & &
     \ar[ll]^{p_S^D}K_D(Y)\ar@{^{(}->}[dd]^{i_{D}^*} \\
     \\
     F(\mathcal{E};R[S]) & & 
   \ar[ll]^{\tilde{\gamma}}F(\mathcal{E};R[D]) } $$
o{\`u} les applications $\tilde{\gamma} :  F(\mathcal{E};
R[D]) \rightarrow 
 F(\mathcal{E};R[S])$ et $p_S^D : K_D(Y) \rightarrow   K_{S}(\Gamma)$
     sont  d{\'e}duites de $\gamma : S \rightarrow D $.

\smallskip

Pour tout couple $(\epsilon, \epsilon') \in \mathcal{E}^2$, on a donc 
$$i_{T}^*(\hat{\mu}^{T}_{\epsilon})(\epsilon')=
i_S^*(\hat{\mu}^{S}_{\epsilon})(\epsilon')
= i_S^* p_S^D({\hat{\mu}^{D}}_{\epsilon})(\epsilon')=
\tilde{\gamma}(\mu^{D}_{\epsilon})(\epsilon').$$

D'apr{\`e}s l'expression de $\mu^{D}_{\epsilon}$ (voir th{\'e}or{\`e}me~\ref{restrictionkb}), 
il suffit de prouver que pour tout $\epsilon \in \mathcal{E}$
et tout $i \in \{ 1, 2, \ldots, N \}$,
\begin{equation} \label{alphalambda}
\alpha_{i}(\epsilon)=-\tau(\lambda_{i}(\epsilon)).
\end{equation}

Cette formule est d{\'e}montr{\'e}e dans \cite{mw4}.

\end{proof}

\subsection{Structure multiplicative}

Comme dans le paragraphe pr{\'e}c{\'e}dent, les r{\'e}sultats de la
section~\ref{multiplicationtours} nous permettent de d{\'e}duire les
th{\'e}or{\`e}mes suivants : 

\begin{theo}\label{smkbs}

La $K$-th{\'e}orie $T$-{\'e}quivariante $K_T(\Gamma)$ de la vari{\'e}t{\'e} de
Bott-Samelson $\Gamma$ s'identifie {\`a} la $R[T]$-alg{\`e}bre 
 $$  R[T][X_1, \ldots , X_N, X_1^{-1},
 \ldots ,  X_N^{-1}]  / \mathcal{J} ,$$ 
o{\`u}  $\mathcal{J}$ est  l'id{\'e}al engendr{\'e} par
$\{ X_i^2-X_i+(1-X_i)M_i \}_{1
  \leq i \leq N}$, avec 
$$\boxed{M_i=e^{-\mu_i}\prod_{j<i}X_j^{-b_{j,i}}.}$$

 \end{theo}

\begin{exem}

On se place dans le cas $A_2$ o{\`u} le groupe $G$ est isomorphe {\`a}
$SL(3,\mathbb{C})$. 
L'alg{\`e}bre $K_{T}\big(\Gamma(\alpha_1,\alpha_2 ,
\alpha_1)\big)$ s'identifie alors {\`a} 
 $$  R[T][X_1, X_2 , X_3, X_1^{-1},
 X_2^{-1} ,  X_3^{-1}]  / \mathcal{J} ,$$ 
o{\`u}  $\mathcal{J}$ est  l'id{\'e}al engendr{\'e} par
$$\{ X_1^2-X_1+(1-X_1)e^{-\alpha_1},  X_2^2-X_2+(1-X_2)e^{-\alpha_2}X_1,
  X_3^2-X_3+(1-X_3)e^{-\alpha_1}X_1^{-2}X_2  \}.$$

 \end{exem}

\medskip

On introduit la $R[T]$-alg{\`e}bre 
$$ \mathcal{T}=R[T][ X_1, \ldots ,  X_N, 
X_1^{-1}, \ldots  , X_N^{-1}, Z_1, \ldots , Z_N].
$$

On note $M=(\mu_1, \ldots , \mu_N)$ le $N$-uplet de racines simples
qui d{\'e}finit la vari{\'e}t{\'e} $\Gamma$.

\begin{defi}

Soit $\epsilon$ un {\'e}l{\'e}ment de  $\mathcal{E}$ de longueur $l$ strictement positive. 
On note $\{i_1 < \cdots <i_l \}$ les {\'e}l{\'e}ments
 de $\pi_+(\epsilon)$. On d{\'e}finit alors l'application $R_M^{\epsilon} : 
\mathcal{T} \rightarrow R[T]$ de la mani{\`e}re suivante :

\begin{enumerate}

\item[$(i)$] $R_M^{\epsilon}$ est $R[T]$-lin{\'e}aire,


\item[$(ii)$] si $P\in \mathcal{T}$ est un mon{\^o}me non nul qui s'{\'e}crit 
sous la forme $P=SX_{i_l}^r Z_{i_l}^s$
  o{\`u} $S \in \mathcal{T}$ est de degr{\'e} $0$ en  $X_{i_l}$
   et $Z_{i_l}$, et o{\`u} $s$ est un entier positif et $r$
  un entier quelconque, 
$$R_M^{\epsilon}(P)=\left\{ \begin{array}{ll}
R_M^{\epsilon-(i_l)}\Big(S(1-M_{i_l})^{s-1}M_{i_l}^r\Big) & { \rm si \,
} s>0, \\
-R_M^{\epsilon-(i_l)}\Big(S(M_{i_l}+M_{i_l}^2 + \cdots + M_{i_l}^{r-1})\Big)
& { \rm si \,} s=0 { \rm  \, \, et \, } r>1, \\
R_M^{\epsilon-(i_l)}\Big(S(1+M_{i_l}^{-1} + \cdots
+ M_{i_l}^{r})\Big)
& { \rm si \,} s=0 { \rm \, \,  et \, } r < 0, \\
R_M^{\epsilon-(i_l)}(S) & { \rm si \,} s=0 { \rm \, \,
  et \, } r = 0, \\
0  & { \rm si \,} s=0 { \rm \, \,
  et \, } r = 1,
\end{array} \right.
$$
\item[$(iii)$] $R_M^{(\mathbf{0})}(P)=P(X_i=1, Z_i=0)$.

\end{enumerate}

\smallskip

Ces trois relations d{\'e}finissent compl{\`e}tement (r{\'e}cursivement) les applications
 $R_M^{\epsilon}$. On pose $R_M=R_M^{(\mathbf{1})}$.

\end{defi}

\begin{theo} \label{cskbs}

Pour $\epsilon \in \mathcal{E}$, on pose $S_{\epsilon}^T=\prod_{i \in
   \pi_-(\epsilon)} X_i \prod_{j \in \pi_+(\epsilon)}Z_j \in
   \mathcal{T}$.
Soient $\epsilon$, $\epsilon'$ et $\epsilon''$ trois {\'e}l{\'e}ments de
$\mathcal{E}$, alors : 
 $$\chi(\overline{\Gamma}_{\epsilon''}, \hat{\mu}^T_{\epsilon} 
 \hat{\mu}^T_{\epsilon'}) =R_M^{\epsilon''}
(S_{\epsilon}^T S_{\epsilon'}^T),$$
 et on a donc : 
$$ \hat{\mu}^T_{\epsilon}  \hat{\mu}^T_{\epsilon'}=
\sum_{\tilde{\epsilon} \in \mathcal{E}} R_M^{\tilde{\epsilon}}
(S_{\epsilon}^T S_{\epsilon'}^T)  \hat{\mu}^T_{\tilde{\epsilon}}.
$$

\end{theo}

\begin{exem}

On se place dans le cas $A_2$, et on prend $M=(\alpha_1, \alpha_2,
\alpha_1)$, $\epsilon=(1,0,0)$, $\epsilon'=(0,0,1)$ et
$\epsilon''=(\mathbf{1})$. On calcule $\chi(\Gamma,  \hat{\mu}^T_{\epsilon} 
 \hat{\mu}^T_{\epsilon'})$ :
 $$\chi(\Gamma, \hat{\mu}^T_{\epsilon} 
 \hat{\mu}^T_{\epsilon'}) =R_M\big( (Z_1X_2X_3)(X_1X_2Z_3 )\big)=
R_M (X_1Z_1X_2^2X_3Z_3 )
$$
$$=R_M^{(\overline{2})}(X_1Z_1X_2^2M_3)=R_M^{(\overline{2})}(X_1Z_1X_2^2
e^{-\alpha_1}X_1^{-2}X_2)=e^{-\alpha_1}R_M^{(\overline{2})}
(X_1^{-1}Z_1X_2^3)
$$
$$=-e^{-\alpha_1}R_M^{(\overline{1})}\big( X_1^{-1}Z_1(M_2+M_2^2)
\big)=
-e^{-\alpha_1}\big[R_M^{(\overline{1})}(X_1^{-1}Z_1e^{-\alpha_2}X_1)-
R_M^{(\overline{1})}(X_1^{-1}Z_1e^{-2\alpha_2}X_1^2)\big]
$$
$$=-e^{-\alpha_1}[e^{-\alpha_2}R_M^{(\overline{1})}(Z_1)-
e^{-2\alpha_2}R_M^{(\overline{1})}(X_1 Z_1)\big]
=-e^{-\alpha_1}[e^{-\alpha_2}R_M^{(\mathbf{0})}(1)-
e^{-2\alpha_2}R_M^{(\mathbf{0})}(e^{-\alpha_1})\big]
$$
$$=-e^{-\alpha_1}(e^{-\alpha_2}-
e^{-2\alpha_2}e^{-\alpha_1})=e^{-2\alpha_1 -  2\alpha_2
}-e^{-\alpha_1-\alpha_2}.
$$
\end{exem}

\section{K-th{\'e}orie {\'e}quivariante des vari{\'e}t{\'e}s de drapeaux} \label{sectionvd}

\subsection{D{\'e}finitions}

La vari{\'e}t{\'e} de drapeaux $X$ n'{\'e}tant pas compacte en g{\'e}n{\'e}ral, on d{\'e}finit $K_T(X)$
de la mani{\`e}re suivante : 

On d{\'e}finit, pour tout entier $n \geq 0$, 
$$\displaystyle{X_{n}=\!\!\!\!\!
 \bigcup_{ \tiny 
\begin{array}{ll} \hspace{0,15 cm} w \in W \\ 
 l(w) \leq n
\end{array}}\!\!\!\!\!X_{w}}.$$

 Soit $\mathcal{F}$ la filtration : 
$$\mathcal{F} :
\emptyset=X_{-1} \subset X_{0} \subset X_{1} \subset \cdots . $$ 

Alors : 

\begin{enumerate}

\item[$(1)$] chaque $X_{n}$ est un sous espace compact $T$-stable de $X$ et,

\item[$(2)$] la topologie de $X$ est la topologie limite induite par la filtration
    $\mathcal{F}$.

\end{enumerate}

\smallskip

Gr{\^a}ce {\`a} cette filtration, on d{\'e}finit alors la K-th{\'e}orie $T$-{\'e}quivariante de $X$, 
not{\'e}e $K_{T}(X)$, par
$$\displaystyle{K_{T}(X)=\lim_{  \leftarrow}}_{ 
 \tiny n \rightarrow +\infty
}\!\!\!\!\! \!\!\!\!\!\!\!\! K_{T}(X_{n}).$$

 Cette d{\'e}finition est ind{\'e}pendante  de la filtration $\mathcal{F}$
v{\'e}rifiant $(1)$ et $(2)$.

\bigskip

On note $F(W; Q[T])$) la 
$R[T]$-alg{\`e}bre des fonctions de $W$ {\`a}
valeurs dans $Q[T]$ munie de
l'addition et de la multiplication point par point. Pour tout $1 \leq i \leq r$,
on d{\'e}finit alors un op{\'e}rateur de Demazure $D_{i}$ sur $F(W; Q[T])$ par 

$$(D_{i}f)(v)=\frac{f(v)-f(vs_{i})e^{-v\alpha_{i}}}{1-e^{-v
\alpha_{i}}}.$$
   
Dans \cite{kkk}, Kostant et Kumar montrent que ces op{\'e}rateurs de Demazure 
v{\'e}rifient les relations de tresses de $W$. Pour tout $w \in W$,
on peut donc d{\'e}finir un op{\'e}rateur $D_{w}$ sur $F(W; Q[T])$  par
$D_w=D_{i_1}D_{i_2} \cdots D_{i_l}$ si $w=s_{i_1}s_{i_2} \cdots s_{i_l}$ est une
d{\'e}composition r{\'e}duite de $w$. 

De plus, pour tout 
$1 \leq i \leq r$, $D_i^2=D_i$. Donc, si pour $\underline{u} \in 
\underline{W}$, on note $D_{\underline{u}}=D_u$, alors pour tout couple 
$(\underline{v},\underline{w})\in \underline{W}^2$, $
D_{\underline{v}}D_{\underline{w}} =D_{\underline{v}\, \underline{w}}$.

\medskip

On note $\Psi$ la
sous-alg{\`e}bre de $F(W; R[T])$ d{\'e}finie par  
$$ \Psi = \{ f \in F(W; R[T]), \hspace{0,2 cm} 
{\rm telles \hspace{0,2 cm} que } \hspace{0,2 cm}
\forall w \in W, \hspace{0,2 cm} D_{w}f \in F(W; R[T]) \}.$$

\medskip

  L'ensemble des points fixes $X^T \approx  W$ 
{\'e}tant discret, on
peut identifier $K_{T}(X^T)$ avec $F(W; R[T])$ et on obtient ainsi un morphisme 
$i_{T}^* : K_{T}(X) \rightarrow F(W; R[T])$. Le r{\'e}sultat suivant est
prouv{\'e} dans \cite{kkk} :

\begin{prop} \label{propositionKKK}

L'application $i_{T}^*$ est injective, et l'image de $K_{T}(X)$ par cette
application est {\'e}gale {\`a} $\Psi$. De plus, $\Psi=\prod_{w \in W}R[T]\psi^w$, o{\`u}
les fonctions $\psi^w$ sont uniquement d{\'e}termin{\'e}es par les relations 
$$\forall (v,w) \in W^2, \, D_v(\psi^w)(1)=\delta_{v,w}.$$

De plus, les fonctions  $\psi^w$
v{\'e}rifient les propri{\'e}t{\'e}s suivantes :  

\smallskip

\begin{enumerate}

 \item[$(i)$] $\psi^{w}(v)=0$ sauf si $w \leq v$,
   
  \item[$(ii)$] $\psi^{w}(w) = \prod_{\beta \in \Delta(w^{-1})}
(1-e^{\beta})$,

 \item[$(iii)$]  $ D_{i}\psi^{w} = \left\{ \begin{array}{ll} \psi^{w}
 + \psi^{ws_{i}} 
 & {\rm si } \hspace{0,2 cm} ws_{i}<w, \\ 
0 & {\rm si } \hspace{0,2 cm} ws_{i}>w,  \end{array} \right.$

\item[$(iv)$] $ \forall v \in W, \psi^{1}(v) =e^{\rho-v\rho} $.
    
\end{enumerate}

\end{prop}

\begin{rema}

Un {\'e}l{\'e}ment $f=(a_{w})_{w \in W}$ de $\prod_{w \in W}R[T]\psi^{w}$ est bien une 
fonction de $W$ {\`a} valeurs dans $R[T]$. En effet soit $v \in W$, d'apr{\`e}s la
propri{\'e}t{\'e}~$(i)$, 
$\sum_{w \in W}a_{w}\psi^{w}(v)$ est une somme finie o{\`u} les termes
{\'e}ventuellement non nuls correspondent aux {\'e}l{\'e}ments $u$ de $W$ qui v{\'e}rifient $u
\leq v$.

\end{rema}




On pose $\hat{\psi}^w=(i_{T}^*)^{-1}(\psi^{w})$, et pour $v \in W$ on note
$\hat{D}_{v} : K_{T}(X) \rightarrow K_{T}(X)$ l'application induite par $D_v :
\Psi \rightarrow \Psi$. La famille $\{\hat{\psi}^w\}_{w \in W}$ est une
base du $R[T]$-module $K_T(X)$.

\begin{rema}

Dans le cas fini, $K_T(X)$ s'identifie {\`a} $K^0(H,X)$ (voir \cite{kkk}), et
Kostant et Kumar montrent dans \cite{kkk} que la base $\{\hat{\psi}^w\}_{w\in
  W}$ de $K_T(X) \simeq K^0(H,X)$ est r{\'e}li{\'e}e aux vari{\'e}t{\'e}s de Schubert par les
relations  
$$ \forall (v,w) \in W^2, \chi(\overline{X}_v, *\hat{\psi}^w)=\delta_{v,w}.
$$
\end{rema}

\medskip 
 
 Dans \cite{kkk},  Kostant et Kumar composent $i_{T}^{*}$ avec $\phi : 
 F(W; Q[T]) \rightarrow
 F(W; Q[T])$ d{\'e}finie par 
 $\phi(f)(w)=f(w^{-1})$ pour tout {\'e}lement $f$ de $F(W; Q[T])$ et tout $w \in
 W$. Ils trouvent alors la sous alg{\`e}bre $\Psi'$ (not{\'e}e $\Psi$ dans \cite{kkk}) 
 de $F(W;R[T])$ : 
 $$\Psi' = \{ f \in F(W;R[T]), \hspace{0,2 cm} {\rm telles \hspace{0,2
 cm} que }  \hspace{0,2 cm}
 \forall w \in W, \hspace{0,2 cm} D'_{w}f \in F(W;R[T]) \},$$
 o{\`u} les op{\'e}rateurs $D_{w}'$ sont d{\'e}finis {\`a} partir des op{\'e}rateurs
 $D_{i}'$ donn{\'e}s  par 
 $$(D_{i}'f)(v)=\frac{f(v)-f(s_{i}v)e^{-v^{-1}\alpha_{i}}}{1-e^{-v^{-1}
 \alpha_{i}}}.$$

 Ils consid{\`e}rent les bases $\psi_{w} '$ de $\Psi'$ 
(not{\'e}e   $\psi^w$ dans \cite{kkk}) et $\tau^w$ de $K_T(X)$
 reli{\'e}es aux  bases $\psi^{w}$ et $\hat{\psi}^w$ consid{\'e}r{\'e}es dans
 cet article  par  les relations
 $\psi_{w}'=\phi(\psi^{w^{-1}})$ et $\tau^w
 =\hat{\psi}^{w^{-1}}$. 

\subsection{Lien avec les vari{\'e}t{\'e}s de Bott-Samelson}

Soit $\mu_{1}, \ldots , \mu_{N}$ une suite quelconque de $N$ racines
simples. On pose $\Gamma=\Gamma(\mu_{1}, \ldots, \mu_{N})$ et on d{\'e}finit une
application $g$ de $\Gamma$ dans $X$ par
multiplication : 
$$g([g_{1},\ldots ,g_{N}]) = g_{1}\times
\cdots   \times g_{N} \; \,[B],$$
o{\`u} $\times$ d{\'e}signe la multiplication dans le groupe $G$. 
Cette application est $T$-{\'e}quivariante.

Le th{\'e}or{\`e}me suivant d{\'e}montr{\'e} dans \cite{mw3} fait le lien entre
$K_T(X)$ et $K_T(\Gamma)$ : 

\begin{theo}  \label{g*k}

Pour tout {\'e}l{\'e}ment $v$ du groupe de Weyl $W$,
$$g^*(\hat{\psi}^{v})=\!\!\!\!\!\!\!\sum_{\epsilon \in \mathcal{E},
 \, \underline{v}(\epsilon)=\underline{v}}\!\!\!\!\!\!\!*
\hat{\mu}_{\epsilon}^T \hspace{0,1 cm}.
$$

\end{theo}

Pour d{\'e}montrer ce th{\'e}or{\`e}me, on a utilis{\'e} dans \cite{mw3} le r{\'e}sultat
suivant dont on va donner une nouvelle d{\'e}monstration : 

\begin{lemm} \label{*chi}

Soit $\hat{\psi} \in K_T(X)$ et $\psi = i_T^*(\hat{\psi})$. Pour tout
$\epsilon \in \mathcal{E}$, on a alors 
$$ \chi(\overline{\Gamma}_{\epsilon},g^*(*\hat{\psi}))=
 *(D_{\underline{v}(\epsilon)}(\psi)(1)).$$

\end{lemm}

\begin{proof}

Dans \cite{mw3}, on a d{\'e}montr{\'e} ce r{\'e}sultat en  g{\'e}n{\'e}ralisant  un
argument de g{\'e}om{\'e}trie alg{\'e}brique
utilis{\'e} par Kostant et Kumar dans \cite{kkk}. On va ici donner une
d{\'e}monstration bas{\'e}e sur un calcul du caract{\`e}re {\`a} l'aide de la formule de
localisation.

\smallskip

 Soit $\hat{\psi} \in K_{T}(X)$, on proc{\`e}de par r{\'e}currence sur $l(\epsilon)$. 
Le r{\'e}sultat est trivial si $l(\epsilon)=0$. 

Supposons le r{\'e}sultat v{\'e}rifi{\'e} pour
tout $\epsilon'$ de longueur strictement inf{\'e}rieure {\`a} $p$, et soit $\epsilon$ de
longueur $p$. D'apr{\`e}s le lemme~\ref{pointsfixestours} et la
formule~\ref{alphalambda}, pour $\epsilon' \leq \epsilon$, les poids
de la repr{\'e}sentation de $H$ dans l'espace tangent {\`a}
$\overline{\Gamma}_{\epsilon}$ en $\epsilon'$ sont les 
$\{-\alpha_{i}(\epsilon') \}_{i\in \pi_+(\epsilon)}$.
En utilisant  la proposition~\ref{pointsfixesab}, on obtient alors 

$$ \chi(
\overline{\Gamma}_{\epsilon},g^*(*\hat{\psi}))= \sum_{\epsilon' \leq
\epsilon}\frac{i_{T}^*(g^*(*\hat{\psi}  ))(\epsilon')}{\prod_{i \in
  \pi_{+}(\epsilon)}
(1-e^{\alpha_{i}(\epsilon')})}= \sum_{\epsilon' \leq
\epsilon}\frac{*\psi(v(\epsilon'))}{\prod_{i \in \pi_{+}(\epsilon)}
(1-e^{\alpha_{i}(\epsilon')})}.$$

Soit $j$ le plus grand {\'e}l{\'e}ment de $\pi_+(\epsilon)$, et soit $\tilde{\epsilon}=
\epsilon -(j)$. En distinguant les {\'e}l{\'e}ments $\epsilon'$ tels que $\epsilon'_j=0$ et
ceux tels que $\epsilon'_j=1$, on obtient 

$  \chi(
\overline{\Gamma}_{\epsilon},g^*(*\hat{\psi} ))=
 \sum_{\epsilon' \leq
\tilde{\epsilon}} \frac{*\psi(v(\epsilon'))}
{(1-e^{\alpha_{j}(\epsilon')})\prod_{i \in
  \pi_{+}(\tilde{\epsilon})}
(1-e^{\alpha_{i}(\epsilon')})}$

\hspace{4 cm} $ + \sum_{\epsilon' \leq
\tilde{\epsilon}}\frac{*\psi(v(\epsilon')s_{\mu_j})}
{(1-e^{-\alpha_{j}(\epsilon')})\prod_{i \in
  \pi_{+}(\tilde{\epsilon})}
(1-e^{\alpha_{i}(\epsilon')})}$.

On a donc  
$$ \chi(
\overline{\Gamma}_{\epsilon},g^*(*\hat{\psi} ))= \sum_{\epsilon' \leq
\tilde{\epsilon}}\frac{1}{\prod_{i \in
  \pi_{+}(\tilde{\epsilon})}
(1-e^{\alpha_{i}(\epsilon')})}
\left[\frac{\psi(v(\epsilon'))-e^{-v(\epsilon')\alpha_{j}}
\psi(v(\epsilon')s_{\mu_j})}
{1-e^{-v(\epsilon')\alpha_{j}}}\right]^*$$
$$= \sum_{\epsilon' \leq
\tilde{\epsilon}}\frac{*i_{T}^*(\hat{D}_{s_{\mu_j}}\hat{\psi} )
(v(\epsilon'))}{\prod_{i \in \pi_{+}(\epsilon)}
(1-e^{\alpha_{i}(\epsilon')})},
$$
i.e., d'apr{\`e}s la formule initiale,
$$ \chi(
\overline{\Gamma}_{\epsilon},g^*(* \hat{\psi} ))=\chi(
\overline{\Gamma}_{\tilde{\epsilon}},g^*(*\hat{D}_{s_{\mu_j}}\hat{\psi} )),$$
et donc par hypoth{\`e}se de r{\'e}currence ($\tilde{\epsilon}$ {\'e}tant de longueur $p-1$),
$$\chi(\overline{\Gamma}_{\epsilon},g^*(* \hat{\psi}  ))=
*(D_{\underline{v}(\tilde{\epsilon})}
i_{T}^*(\hat{D}_{s_{\mu_j}}\hat{\psi}  ))(1)=*(D_{\underline{v}(\tilde{\epsilon})}
D_{s_{\mu_j}}(\psi ))(1)=*(D_{\underline{v}(\epsilon)}(\psi))(1).$$

\end{proof}

\begin{proof}[D{\'e}monstration du th{\'e}or{\`e}me~\ref{g*k}]

On termine la d{\'e}monstration du th{\'e}or{\`e}me de la m{\^e}me mani{\`e}re que dans \cite{mw3}.

Soit $v$ un {\'e}l{\'e}ment du groupe de Weyl $W$.
D'apr{\`e}s le lemme~\ref{*chi}, pour tout {\'e}l{\'e}ment $\epsilon \in
\mathcal{E}$,
$$\chi(
\overline{\Gamma}_{\epsilon},g^*(*\hat{\psi}^v))=
*(D_{\underline{v}(\epsilon)})(\psi^v)(1).$$

 Or, d'apr{\`e}s la caract{\'e}risation des fonctions 
$\{\psi^w\}_{w \in W}$ (proposition~\ref{propositionKKK}),
$$\forall u \in W, (D_{u}(\psi^{v}))(1)=\delta_{u,v}.$$

On d{\'e}duit des deux formules pr{\'e}c{\'e}dentes que pour tout $\epsilon \in
\mathcal{E}$,
\begin{equation} \label{chi}
  \chi(\overline{\Gamma}_{\epsilon},g^*(*\hat{\psi}^v))=
\delta_{\underline{v}(\epsilon)
 ,\underline{v}}. \hspace{2 cm} \end{equation}

D'apr{\`e}s la caract{\'e}risation de la base
$\{\hat{\mu}_{\epsilon}^T\}_{\epsilon \in \mathcal{E}}$, on a donc
bien 
$$g^*(\hat{\psi}^{v})=\!\!\!\!\!\!\!\sum_{\epsilon \in \mathcal{E},
 \, \underline{v}(\epsilon)=\underline{v}}\!\!\!\!\!\!\!*
\hat{\mu}_{\epsilon}^T \hspace{0,1 cm}.$$

\end{proof}

\begin{rema}

Ce th{\'e}or{\`e}me est utilis{\'e} dans \cite{mw3} pour donner une formule
explicite pour les fonctions $\psi^w$. Cette formule est d{\'e}montr{\'e}e
par d'autres m{\'e}thodes par William Graham dans \cite{gra}.

\end{rema}

\begin{exem}

On se place dans le cas $A_2$, et on consid{\`e}re
$\Gamma=~\Gamma(\alpha_1, \alpha_2, \alpha_1)$. Le
th{\'e}or{\`e}me~\ref{g*k} nous donne les relations suivantes :

$$
\begin{array}{lcl}
g^*(\hat{\psi}^1) & =  & *\hat{\mu}_{(\bf{0})}^T, \\
g^*(\hat{\psi}^{s_1}) & = &
*\hat{\mu}_{(1,0,0)}^T+*\hat{\mu}_{(0,0,1)}^T+*\hat{\mu}_{(1,0,1)}^T, \\
g^*(\hat{\psi}^{s_2}) & = & *\hat{\mu}_{(0,1,0)}^T, \\
g^*(\hat{\psi}^{s_1s_2})&  = & *\hat{\mu}_{(1,1,0)}^T, \\
g^*(\hat{\psi}^{s_2s_1}) & = & *\hat{\mu}_{(0,1,1)}^T, \\
g^*(\hat{\psi}^{s_1s_2s_1}) & = & *\hat{\mu}_{(\bf{1})}^T.
\end{array}
$$

\end{exem}

\subsection{Structure multiplicative}

On note $q_{u,v}^w \in R[T]$ les constantes de structure d{\'e}finies par
les relations  
\begin{equation} \label{constantesschubert}
 \hat{\psi}^u \hat{\psi}^v = \sum_{ w \in W}
 q_{u,v}^w \hat{\psi}^w.
\end{equation}

\begin{rema} \label{nullite}

M{\^e}me dans le cas infini, cette somme a bien un sens. En effet, pour tout
$u' \in W$, les classes $\hat{\psi}^w$ restreintes {\`a}
$\overline{X}_{u'}$ sont nulles sauf si $ w\leq u'$.

\end{rema}

Dans \cite{kkk}, Kostant et Kumar donnent une formule pour calculer
ces constantes de structure. Ils d{\'e}finissent la matrice
$E=(e^{u,v})_{(u,v)\in W^2}$ par $$e^{u,v}=\psi^u(v),$$ et pour $w \in
W$ la matrice diagonale $E_w=(E_w(u,v))_{(u,v)\in W^2}$ par
$$E_w(u,v)=\delta_{u,v}\psi^w(u).$$

Pour $w\in W$, on note $Q_{w}=(Q_w(u,v))_{(u,v)\in W^2}$ la matrice
d{\'e}finie par $$Q_w(u,v)=q_{u,v}^w.$$

 Kostant et Kumar montrent alors la formule suivante :
$$Q_w=EE_wE^{-1},
$$
o{\`u} l'inverse $E^{-1}$ de $E$ est une matrice {\`a} coefficients dans
$Q[T]$. Pour trouver un coefficient, on doit donc calculer la matrice
$E$ et son inverse. On obtient ainsi des expressions dans $Q[T]$ qu'il
faut simplifier puisque la matrice $Q_w$ est {\`a} coefficients dans
$R[T]$. On va donner une m{\'e}thode de calcul plus efficace si on veut
calculer un coefficient particulier, et qui ne passe pas par le corps
des fractions de $R[T]$.

\bigskip

On fixe une d{\'e}composition r{\'e}duite $w=s_{\mu_1} \cdots s_{\mu_N}$
d'un {\'e}l{\'e}ment $w$ du groupe de Weyl $W$. On pose $M=(\mu_1, \ldots,
\mu_N)$.

\begin{theo} \label{cskvd}

Pour tout couple $(u,v) \in W^2$, 
$$q_{u,v}^w = *R_M\Big((\!\!\!\!\sum_{\epsilon \in \mathcal{E},
 \,
 \underline{v}(\epsilon)=\underline{u}}\!\!\!\!\!\!\!S_{\epsilon}^T)
 \,\,(\!\!\!\!\sum_{\epsilon' \in \mathcal{E},
 \,
 \underline{v}(\epsilon')=\underline{v}}\!\!\!\!\!\!\!S_{\epsilon'}^T) 
\Big).
$$

\end{theo}

\begin{proof}

On pose $\Gamma=\Gamma(\mu_{1}, \ldots, \mu_{N})$, et on rappelle la
d{\'e}finition  de l'application
 $g$ de $\Gamma$ dans $X$ : 
$$g([g_{1},\ldots ,g_{N}]) = g_{1}\times
\cdots   \times g_{N} \; \,[B].$$

Si on applique $g^*$ {\`a}
l'{\'e}galit{\'e}~\ref{constantesschubert}, on obtient  
$$  g^*(\hat{\psi}^u\hat{\psi}^v) = \sum_{ \tilde{w} \in W}
 q_{u,v}^{\tilde{w}} g^*  (\hat{\psi}^{\tilde{w} }),
$$
d'o{\`u}, d'apr{\`e}s le th{\'e}or{\`e}me~\ref{g*k},
 $$ ( \!\!\!\!\sum_{\epsilon \in \mathcal{E},
 \, \underline{v}(\epsilon)=\underline{u}}\!\!\!\!\!\!\!*
\hat{\mu}_{\epsilon}^T)(\!\!\!\!\sum_{\epsilon' \in \mathcal{E},
 \, \underline{v}(\epsilon')=\underline{v}}\!\!\!\!\!\!\!*
\hat{\mu}_{\epsilon'}^T  )  = \sum_{ \tilde{w} \in W}
 q_{u,v}^{\tilde{w}}  \!\!\!\!\sum_{\tilde{\epsilon} \in \mathcal{E},
 \, \underline{v}(\tilde{\epsilon})=\underline{\tilde{w}}}\!\!\!\!\!\!\!*
\hat{\mu}_{\tilde{\epsilon}}^T,
$$
d'o{\`u} l'on tire
 $$ ( \!\!\!\!\sum_{\epsilon \in \mathcal{E},
 \, \underline{v}(\epsilon)=\underline{u}}\!\!\!\!\!\!\!
\hat{\mu}_{\epsilon}^T)(\!\!\!\!\sum_{\epsilon' \in \mathcal{E},
 \, \underline{v}(\epsilon')=\underline{v}}\!\!\!\!\!\!\!
\hat{\mu}_{\epsilon'}^T  )  = \sum_{ \tilde{w} \in W}
*q_{u,v}^{\tilde{w}}  \!\!\!\!\sum_{\tilde{\epsilon} \in \mathcal{E},
 \, \underline{v}(\tilde{\epsilon})=\underline{\tilde{w}}}\!\!\!\!\!\!\!
\hat{\mu}_{\tilde{\epsilon}}^T.
$$

Dans le terme de droite, le coefficient de $\hat{\mu}_{(\mathbf1)}^T $
est {\'e}gal {\`a} $$*q_{u,v}^{w}.$$ 

Dans le terme de gauche, d'apr{\`e}s le
th{\'e}or{\`e}me~\ref{cskbs},  le coefficient de $\hat{\mu}_{(\mathbf1)}^T $
est {\'e}gal {\`a} 
$$R_M\Big((\!\!\!\!\sum_{\epsilon \in \mathcal{E},
 \,
 \underline{v}(\epsilon)=\underline{u}}\!\!\!\!\!\!\!S_{\epsilon}^T)
 \,\,(\!\!\!\!\sum_{\epsilon' \in \mathcal{E},
 \,
 \underline{v}(\epsilon')=\underline{v}}\!\!\!\!\!\!\!S_{\epsilon'}^T) 
\Big). $$

La famille $\{ \hat{\mu}_{\tilde{\epsilon}}^T \}_{\tilde{\epsilon} \in
\mathcal{E}}$
{\'e}tant une base du $R[T]$-module $K_T(\Gamma)$, on en d{\'e}duit l'{\'e}galit{\'e} :
 $$ q_{u,v}^{w} = *R_M\Big((\!\!\!\!\sum_{\epsilon \in \mathcal{E},
  \,
  \underline{v}(\epsilon)=\underline{u}}\!\!\!\!\!\!\!S_{\epsilon}^T)
  \,\,(\!\!\!\!\sum_{\epsilon' \in \mathcal{E},
  \,
  \underline{v}(\epsilon')=\underline{v}}\!\!\!\!\!\!\!S_{\epsilon'}^T) 
 \Big).
$$
 \end{proof}

\section{Exemples}  \label{sectionex}

\subsection{Le cas $SL(3,\mathbb{C})$}

On se place ici dans le cas $A_2$, o{\`u} le groupe $G$ est isomorphe {\`a}
$SL(3,\mathbb{C})$ et $B \subset G$ au sous-groupe de
$SL(3,\mathbb{C})$ form{\'e} des matrices triangulaires sup{\'e}rieures. 
Le groupe de Weyl $W$ s'identifie au groupe des
permutations  de l'ensemble $\{1,2,3\}$.

\subsubsection{Calcul du produit $\hat{\psi}^1\hat{\psi}^1$}

Soit $w=s_{\mu_1} \cdots s_{\mu_N}$ une d{\'e}composition
r{\'e}duite d'un {\'e}l{\'e}ment du groupe de Weyl $W$, d'apr{\`e}s le
th{\'e}or{\`e}me~\ref{cskvd} et quelque soit le groupe $G$,
$$q_{1,1}^w=*R_M(X_1^2 X_2^2 \cdots X_N^2),
$$
o{\`u} $M=(\mu_1, \mu_2, \ldots \mu_N)$. En effet, le seul {\'e}l{\'e}ment
$\epsilon$ de $\mathcal{E}$ tel que $\underline{v}(\epsilon)=1$ est
$\epsilon=(\bf{0})$. 

\medskip

Tout d'abord $q_{1,1}^1=1$ (on a toujours $q_{w,w}^w=\psi^w(w)$).

\smallskip

Ensuite,
$$q_{1,1}^{s_1}=*R_{\alpha_1}(X_1^2)=-*R_{\alpha_1}^{(\bf{0})}
(e^{-\alpha_1})=-e^{\alpha_1},$$ 
et de m{\^e}me, $q_{1,1}^{s_2}=-e^{\alpha_2}$.

\smallskip

Ensuite, $$q_{1,1}^{s_1s_2}=*R_{\alpha_1, \alpha_2}(X_1^2
X_2^2)=-*R_{\alpha_1,\alpha_2}^{(\overline{1})}(X_1^2
e^{-\alpha_2}X_1)=-e^{\alpha_2}*R_{\alpha_1, \alpha_2}^{(\overline{1})}(X_1^3)
$$
$$=e^{\alpha_2}*R_{\alpha_1,
  \alpha_2}^{(\bf{0})}(e^{-\alpha_1}+(e^{-\alpha_1})^2)
=e^{\alpha_2}(e^{\alpha_1}+e^{2\alpha_1})=e^{\alpha_1
  +\alpha_2}(1+e^{\alpha_1}),$$ 
et de m{\^e}me,
  $q_{1,1}^{s_2s_1}=e^{\alpha_1+\alpha_2}(1+e^{\alpha_2})$.

\smallskip

Enfin, $q_{1,1}^{s_1s_2s_1}=*R_M(X_1^2 X_2^2 X_3^2)$, o{\`u}
$M=(\alpha_1, \alpha_2, \alpha_1)$. On obtient donc 
$$q_{1,1}^{s_1s_2s_1}=-*R_{M}^{(\overline{2})}
(X_1^2X_2^2e^{-\alpha_1}X_1^{-2}X_2)= 
-e^{\alpha_1}*R_M^{(\overline{2})}(X_2^3)$$
$$=e^{\alpha_1}*R_M^{(\overline{1})}(e^{-\alpha_2}X_1+
(e^{-\alpha_2}X_1)^2)=-e^{\alpha_1+2\alpha_2}*R_M^{(\bf{0})}(e^{-\alpha_1})
=-e^{2\alpha_1+2\alpha_2}.
$$

 On a ainsi obtenu 
$$(\hat{\psi}^1)^2=\hat{\psi}^1 -e^{\alpha_1}\hat{\psi}^{s_1}
-e^{\alpha_2}\hat{\psi}^{s_2}
+e^{\alpha_1+\alpha_2}(1+e^{\alpha_1})\hat{\psi}^{s_1s_2} $$
$$+e^{\alpha_1+\alpha_2}(1+e^{\alpha_2})
\hat{\psi}^{s_2s_1}-e^{2\alpha_1+2\alpha_2}
\hat{\psi}^{s_1s_2s_1}. 
$$

\subsubsection{Calcul du produit $\hat{\psi}^1 \hat{\psi}^{s_1}$} 

Tout d'abord, $q_{1,s_1}^{s_1}=\psi^{1}(s_1)=e^{\alpha_1}$.

\smallskip

Ensuite, $$q_{1,s_1}^{s_1s_2}=*R_{\alpha_1, \alpha_2}(X_1Z_1
X_2^2)=-*R_{\alpha_1,\alpha_2}^{(\overline{1})}(X_1Z_1
e^{-\alpha_2}X_1)$$
$$=-e^{\alpha_2}*R_{\alpha_1, \alpha_2}^{(\overline{1})}(X_1^2Z_1)
=-e^{\alpha_2}*R_{\alpha_1,
  \alpha_2}^{(\bf{0})}(e^{-2\alpha_1})
=-e^{2\alpha_1+\alpha_2},$$ 
et de plus,
  $$q_{1,s_1}^{s_2s_1}=*R_{\alpha_2, \alpha_1}(X_1^2
X_2Z_2)=*R_{\alpha_2,\alpha_1}^{(\overline{1})}(X_1^2
e^{-\alpha_1}X_1)$$
$$=e^{\alpha_1}*R_{\alpha_2, \alpha_1}^{(\overline{1})}(X_1^3)
=-e^{\alpha_1}*R_{\alpha_2,
  \alpha_1}^{(\bf{0})}(e^{-\alpha_2}+e^{-2\alpha_2})
=-e^{\alpha_1+\alpha_2}(1+e^{\alpha_2}).$$

\smallskip

Enfin, $q_{1,s_1}^{s_2s_1s_2}=*R_M(X_1^2 X_2 Z_2 X_3^2)$, o{\`u}
$M=(\alpha_2, \alpha_1, \alpha_2)$. On a donc
$$q_{1,s_1}^{s_2s_1s_2}=-*R_{M}^{(\overline{2})} (X_1^2X_2Z_2
e^{-\alpha_2}X_1^{-2}X_2)= 
-e^{\alpha_2}*R_M^{(\overline{2})}(X_2^2 Z_2)$$
$$=-e^{\alpha_2}*R_M^{(\overline{1})}(e^{-2\alpha_1}X_1^2)=
  e^{2\alpha_1+\alpha_2}* R_M^{(\bf{0})}(e^{-\alpha_2})
=e^{2\alpha_1+2\alpha_2}.
$$

 On a ainsi obtenu 
$$\hat{\psi}^1 \hat{\psi}^{s_1}= e^{\alpha_1}\hat{\psi}^{s_1}
-e^{2\alpha_1+\alpha_2} \hat{\psi}^{s_1s_2} 
-e^{\alpha_1+\alpha_2}(1+e^{\alpha_2})
\hat{\psi}^{s_2s_1}+e^{2\alpha_1+2\alpha_2}
\hat{\psi}^{s_1s_2s_1}. 
$$

On a de m{\^e}me (en changeant $\alpha_1$ en $\alpha_2$ et $s_1$ en
$s_2$) 
$$\hat{\psi}^1 \hat{\psi}^{s_2}= e^{\alpha_2}\hat{\psi}^{s_2}
-e^{\alpha_1+2\alpha_2} \hat{\psi}^{s_2s_1} 
-e^{\alpha_1+\alpha_2}(1+e^{\alpha_1})
\hat{\psi}^{s_1s_2}+e^{2\alpha_1+2\alpha_2}
\hat{\psi}^{s_1s_2s_1}. 
$$

\subsubsection{Calcul du produit $\hat{\psi}^{s_1} \hat{\psi}^{s_1}$}

Tout d'abord, $q_{s_1,s_1}^{s_1}=\psi^{s_1}(s_1)=1-e^{\alpha_1}$.

\smallskip

Ensuite, $$q_{s_1,s_1}^{s_1s_2}=*R_{\alpha_1, \alpha_2}(Z_1^2
X_2^2)=-*R_{\alpha_1,\alpha_2}^{(\overline{1})}(Z_1^2
e^{-\alpha_2}X_1)
=-e^{\alpha_2}*R_{\alpha_1, \alpha_2}^{(\overline{1})}(X_1Z_1^2)$$
$$=-e^{\alpha_2}*R_{\alpha_1,
  \alpha_2}^{(\bf{0})}(e^{-\alpha_1}(1-e^{-\alpha_1}))
=-e^{\alpha_1+\alpha_2}(1-e^{\alpha_1}),$$ 
et de plus,
  $$q_{s_1,s_1}^{s_2s_1}=*R_{\alpha_2, \alpha_1}(X_1^2
Z_2^2)=*R_{\alpha_2,\alpha_1}^{(\overline{1})}(X_1^2
(1-e^{-\alpha_1}X_1))
=*R_{\alpha_2, \alpha_1}^{(\overline{1})}(X_1^2-e^{-\alpha_1}X_1^3)$$
$$=-*R_{\alpha_2,
  \alpha_1}^{(\bf{0})}(e^{-\alpha_2}-e^{-\alpha_1}(e^{-\alpha_2} 
+e^{-2\alpha_2 }))
=-e^{\alpha_2}(1-e^{\alpha_1}-e^{\alpha_1+\alpha_2}).$$

\smallskip

Enfin, $q_{s_1,s_1}^{s_2s_1s_2}=*R_M(X_1^2  Z_2^2 X_3^2)$, o{\`u}
$M=(\alpha_2, \alpha_1, \alpha_2)$. On a donc
$$q_{s_1,s_1}^{s_2s_1s_2}=-*R_{M}^{(\overline{2})} (X_1^2Z_2^2
e^{-\alpha_2}X_1^{-2}X_2)= 
-e^{\alpha_2}*R_M^{(\overline{2})}(X_2 Z_2^2)$$
$$=-e^{\alpha_2}*R_M^{(\overline{1})}(e^{-\alpha_1}X_1
( 1-e^{-\alpha_1}X_1 ))=
  e^{\alpha_1+\alpha_2}* R_M^{(\bf{0})}(-e^{-\alpha_1-\alpha_2})
=-e^{2\alpha_1+2\alpha_2}.
$$

  On a ainsi obtenu 
$$(\hat{\psi}^{s_1})^2= (1-e^{\alpha_1})\hat{\psi}^{s_1}
-e^{\alpha_1+\alpha_2}(1-e^{\alpha_1}) \hat{\psi}^{s_1s_2} 
-e^{\alpha_2}(1-e^{\alpha_1}-e^{\alpha_1+\alpha_2})
\hat{\psi}^{s_2s_1}-e^{2\alpha_1+2\alpha_2}
\hat{\psi}^{s_1s_2s_1}. 
$$

On a de m{\^e}me (en changeant $\alpha_1$ en $\alpha_2$ et $s_1$ en
$s_2$) :
$$(\hat{\psi}^{s_2})^2= (1-e^{\alpha_2})\hat{\psi}^{s_2}
-e^{\alpha_1+\alpha_2}(1-e^{\alpha_2}) \hat{\psi}^{s_2s_1} 
-e^{\alpha_1}(1-e^{\alpha_2}-e^{\alpha_1+\alpha_2})
\hat{\psi}^{s_1s_2}-e^{2\alpha_1+2\alpha_2}
\hat{\psi}^{s_1s_2s_1}. 
$$

\subsubsection{Calcul du produit $\hat{\psi}^{s_1} \hat{\psi}^{s_2}$}

Tout d'abord,  $$q_{s_1,s_2}^{s_1s_2}=*R_{\alpha_1, \alpha_2}(X_1Z_1
X_2Z_2)=*R_{\alpha_1,\alpha_2}^{(\overline{1})}(X_1Z_1
e^{-\alpha_2}X_1)$$
$$=e^{\alpha_2}*R_{\alpha_1, \alpha_2}^{(\overline{1})}(X_1^2Z_1)
=e^{\alpha_2}*R_{\alpha_1,
  \alpha_2}^{(\bf{0})}(e^{-2\alpha_1})
=e^{2\alpha_1+\alpha_2},$$ 
et de m{\^e}me,
  $q_{s_1,s_2}^{s_2s_1}=e^{\alpha_1+2\alpha_2}.$

\smallskip

De plus, $q_{s_1,s_2}^{s_2s_1s_2}=*R_M \big(X_1Z_2 X_3(X_1X_2Z_3 +
Z_1X_2X_3 + Z_1X_2Z_3) \big)$, o{\`u}
$M=(\alpha_2, \alpha_1, \alpha_2)$. On voit facilement que les deux
derniers termes se compensent, et on obtient 
$$q_{s_1,s_2}^{s_2s_1s_2}=*R_M(X_1^2X_2Z_2X_3Z_3)=
*R_{M}^{(\overline{2})}(X_1^2X_2Z_2e^{-\alpha_2}X_1^{-2}X_2)$$
$$=e^{\alpha_2} *R_M^{(\overline{2})}(X_2^2 Z_2)
=e^{\alpha_2}*R_M^{(\overline{1})}(e^{-2\alpha_1}X_1^2)=
  -e^{2\alpha_1+\alpha_2}* R_M^{(\bf{0})}(e^{-\alpha_2})
=-e^{2\alpha_1+2\alpha_2}.
$$

    On a ainsi obtenu 
$$\hat{\psi}^{s_1}\hat{\psi}^{s_2}= 
e^{2\alpha_1+\alpha_2} \hat{\psi}^{s_1s_2} 
+e^{\alpha_1+2\alpha_2}
\hat{\psi}^{s_2s_1}-e^{2\alpha_1+2\alpha_2}
\hat{\psi}^{s_1s_2s_1}. 
$$

   \subsubsection{Autres r{\'e}sultats}

     On obtient facilement les r{\'e}sultats suivants : 

$$\hat{\psi}^{1}\hat{\psi}^{s_1s_2s_1}= 
e^{2\alpha_1+2\alpha_2}\hat{\psi}^{s_1s_2s_1},$$
$$ \hat{\psi}^{s_1}\hat{\psi}^{s_1s_2s_1}= 
 e^{\alpha_1+\alpha_2}(1-e^{\alpha_1+\alpha_2})\hat{\psi}^{s_1s_2s_1},$$
$$\hat{\psi}^{s_1s_2}\hat{\psi}^{s_1s_2s_1}= 
e^{\alpha_2}(1-e^{\alpha_1})(1-e^{\alpha_1+\alpha_2})\hat{\psi}^{s_1s_2s_1},$$
$$(\hat{\psi}^{s_1s_2s_1})^2 = 
(1-e^{\alpha_1})(1-e^{\alpha_2})(1-e^{\alpha_1+\alpha_2})\hat{\psi}^{s_1s_2s_1},
$$
$$(\hat{\psi}^{s_1s_2})^2 =
(1-e^{\alpha_1})(1-e^{\alpha_1+\alpha_2}) \hat{\psi}^{s_1s_2} 
-e^{\alpha_2}(1-e^{\alpha_1})(1-e^{\alpha_1+\alpha_2})\hat{\psi}^{s_1s_2s_1},
$$
$$\hat{\psi}^{s_1}\hat{\psi}^{s_2s_1}= 
e^{\alpha_2}(1-e^{\alpha_1+\alpha_2}) \hat{\psi}^{s_2s_1} 
-e^{\alpha_1+\alpha_2}(1- e^{\alpha_1+\alpha_2} )\hat{\psi}^{s_1s_2s_1},
$$
$$\hat{\psi}^{s_1}\hat{\psi}^{s_1s_2}= 
e^{\alpha_1+\alpha_2}(1-e^{\alpha_1}) \hat{\psi}^{s_1s_2} 
+e^{2\alpha_1+2\alpha_2}\hat{\psi}^{s_1s_2s_1},
$$
$$\hat{\psi}^{1}\hat{\psi}^{s_1s_2}= 
e^{2\alpha_1+\alpha_2} \hat{\psi}^{s_1s_2} 
-e^{2\alpha_1+2\alpha_2}
\hat{\psi}^{s_1s_2s_1},$$ 
$$\hat{\psi}^{s_1s_2}\hat{\psi}^{s_2s_1}=
e^{\alpha_1+\alpha_2}(1-e^{\alpha_1+\alpha_2})\hat{\psi}^{s_1s_2s_1}.
$$

Les six autres produits s'obtiennent en permutant $s_1$ et $s_2$ d'une
part, et $\alpha_1$ et $\alpha_2$ d'autre part.

\subsection{Quelques autres calculs}

\subsubsection{Calculs dans le cas $B_2$} On se place dans le cas
$B_2$ o{\`u} la matrice de Cartan est $A=\begin{pmatrix} 2 & -2 \\ -1 & 2 
\end{pmatrix}$.

On calcule d'abord $q_{1,s_1}^{s_2s_1s_2}=*R_M(X_1^2X_2Z_2X_3^2)$,
o{\`u} $M=(\alpha_2, \alpha_1, \alpha_2)$. On a donc  
$$q_{1,s_1}^{s_2s_1s_2}=-*R_M^{(\overline{2})}
(X_1^2X_2Z_2e^{-\alpha_2}X_1^{-2}X_2^2) 
=-e^{\alpha_2}*R_M^{(\overline{2})}(X_2^3Z_2)
$$
$$=-e^{\alpha_2}*R_M^{(\overline{1})}(e^{-3\alpha_1}X_1^3)=-e^{3\alpha_1+\alpha_2}
*R_M^{(\overline{1})}(X_1^3)=e^{3\alpha_1+2\alpha_2}(1+e^{\alpha_2}).
$$

 \medskip

 On calcule maintenant
$$q_{1,s_1}^{s_2s_1s_2s_1}=*R_M\big(X_1X_2X_3X_4(X_1Z_2X_3X_4+ X_1Z_2X_3Z_4
+X_1X_2X_3Z_4)\big),$$ o{\`u} $M=(\alpha_2,\alpha_1, \alpha_2,
\alpha_1)$. On voit facilement que les deux premiers termes se
compensent et on obtient 
$$q_{1,s_1}^{s_2s_1s_2s_1}=*R_M(X_1^2X_2^2X_3^2X_4Z_4)=
*R_M^{(\overline{3})}(X_1^2X_2^2X_3^2e^{-\alpha_1}X_1X_2^{-2}X_3)$$
$$=e^{\alpha_1}*R_M^{(\overline{3})}(X_1^3X_3^3)=-e^{\alpha_1}
*R_M^{(\overline{2})}\big(X_1^3(e^{-\alpha_2}X_1^{-2}X_2^2+
e^{-2\alpha_2}X_1^{-4}X_2^4)\big) 
$$
$$=-e^{\alpha_1+\alpha_2}*R_M^{(\overline{2})}(X_1X_2^2+e^{-\alpha_2}X_1^{-1}X_2^4)$$
$$=e^{\alpha_1+\alpha_2}*R_M^{(\overline{1})}(e^{-\alpha_1}X_1^{2}+
e^{-\alpha_1-\alpha_2}+e^{-2\alpha_1-\alpha_2}
X_1+e^{-3\alpha_1-\alpha_2}X_1^2)
$$
$$=e^{\alpha_1+\alpha_2}(-e^{\alpha_1+\alpha_2}+e^{\alpha_1+ \alpha_2}
-e^{3\alpha_1+2\alpha_2}) =-e^{4\alpha_1+3\alpha_2}.
$$

Dans cet exemple, on s'aper{\c c}oit que des termes peuvent se
compenser m{\^e}me quand on n'a qu'un mon{\^o}me au d{\'e}part.

\medskip

On calcule enfin
$q_{s_1s_2,s_1s_2}^{s_2s_1s_2s_1}=*R_M(X_1^2Z_2^2Z_3^2X_4^2)$, o{\`u}
$M=(\alpha_2,\alpha_1,\alpha_2, \alpha_1)$. On obtient 
$$q_{s_1s_2,s_1s_2}^{s_2s_1s_2s_1}=-*R_M^{(\overline{3})}(X_1^2Z_2^2Z_3^2 
e^{-\alpha_1}X_1X_2^{-2}X_3)=-e^{\alpha_1}*R_M^{(\overline{3})}(X_1^3
X_2^{-2}Z_2^2X_3Z_3^2)
$$
$$=-e^{\alpha_1+\alpha_2}*R_M^{(\overline{2})}(X_1Z_2^2(1-e^{-\alpha_2}X_1^{-2}X_2^2))$$
$$=-e^{\alpha_1+\alpha_2}*R_M^{(\overline{1})}(X_1-e^{-\alpha_1}X_1^2-e^{-2\alpha_1-
 \alpha_2}X_1+e^{-3\alpha_1-\alpha_2}X_1^2)$$
$$=-e^{2\alpha_1+2\alpha_2}
(1-e^{2\alpha_1+\alpha_2}).
$$

\subsubsection{Un calcul dans le cas $G_2$} On se place dans le cas
$G_2$ o{\`u} la matrice de Cartan est $A=\begin{pmatrix} 2 & -1 \\ -3 & 2 
\end{pmatrix}$ et on calcule 
$$q_{s_2,s_2s_1}^{s_1s_2s_1s_2}=*R_M\big(
(X_1Z_2X_3X_4+X_1Z_2X_3Z_4+X_1X_2X_3Z_4) 
X_1Z_2Z_3X_4 \big),$$
o{\`u} $M=(\alpha_1,\alpha_2,\alpha_1,\alpha_2)$. On voit facilement que
les deux derniers termes s'annulent et on obtient 
$$q_{s_2,s_2s_1}^{s_1s_2s_1s_2}=*R_M(X_1^2X_2Z_2X_3Z_3X_4Z_4)=
e^{\alpha_2}*R_M^{(\overline{3})}(X_1^3X_2^{-1}Z_2X_3^2Z_3)
$$
$$=e^{2\alpha_1+\alpha_2}*R_M^{(\overline{2})}(X_1^{-1}X_2^5Z_2)=
e^{2\alpha_1+6\alpha_2}*R_M^{(\overline{1})}(X_1^4)=-e^{3\alpha_1+6\alpha_2} 
(1+e^{\alpha_1}+e^{2\alpha_1}).
$$

\subsection{K-th{\'e}orie ordinaire}

La $K$-th{\'e}orie ordinaire de $X$, not{\'e}e $K(X)$, est le groupe construit
{\`a} partir du semi-groupe des classes d'isomorphisme de fibr{\'e}s
vectoriels complexes de dimension finie au dessus de $X$. La
$K$-th{\'e}orie du point s'identifie {\`a} $\mathbb{Z}$, et $K(X)$ est
munie d'une structure d'anneau d{\'e}finie {\`a} l'aide du produit
tensoriel. On note $ev$ l'application canonique $K_T(X) \rightarrow
K(X)$.
Dans \cite{kkk}, Kostant et Kumar montrent le r{\'e}sultat suivant : 
\begin{prop} \label{ordinaire}

L'application canonique $\hat{ev} : \mathbb{Z}\otimes_{R[T]}K_T(X)
\rightarrow K(X) $ est un isomorphisme, o{\`u} $\mathbb{Z}$ est
consid{\'e}r{\'e} comme un $R[T]$-module par l'application de $R[T]$ dans
$\mathbb{Z}$ d{\'e}finie par l'{\'e}valuation en $1$.

\end{prop}

Pour $w \in W$, on pose $\hat{\psi}^w_1=ev(\hat{\psi}^w) \in K(X)$. D'apr{\`e}s la
proposition pr{\'e}c{\'e}dente, la famille $\{ \hat{\psi}^w_1 \}_{w \in W}$ est une
base du $\mathbb{Z}$-module $K(X)$. C'est la base duale (pour $*$) de
la base de $K(X)$ construite par Demazure dans \cite{demazure}
(voir \cite{kkk}). Le th{\'e}or{\`e}me~\ref{cskvd} et la
proposition~\ref{ordinaire} permettent de calculer les constantes de
structure par rapport {\`a} cette base. On note $t_{u,v}^w$ les
entiers relatifs d{\'e}finis par 
$$ \hat{\psi}^u_1 \hat{\psi}^v_1 =\sum_{w \in W}t_{u,v}^w
\hat{\psi}^w_1 . 
$$

Ces entiers se calculent en {\'e}valuant en $1$ les constantes de
structure $q_{u,v}^w$ de $K_T(X)$. On peut ainsi restreindre le
th{\'e}or{\`e}me~\ref{cskvd} au cas de la $K$-th{\'e}orie ordinaire. Pour
cela, on va introduire des applications $r_M$ qui sont les
{\'e}valuations en $1$ des applications $R_M$ d{\'e}finies dans la
section~\ref{sectionBott}.

Soit $w=s_{\mu_1} \cdots s_{\mu_N}$ une d{\'e}composition r{\'e}duite d'un
{\'e}l{\'e}ment $w$ du groupe de Weyl $W$. On pose $M=(\mu_1, \ldots,
\mu_N)$ et
$$\mathcal{U}=\mathbb{Z}[X_1, \ldots, X_N, X_1^{-1}, \ldots ,X_N^{-1},
Z_1, \ldots ,Z_N].$$ 

Pour $1 \leq i \leq N$, on d{\'e}finit $m_i \in
\mathcal{U}$ par 
$$\boxed{m_i = \prod_{j<i}X_j^{-\mu_i(\mu_j^{\vee})}.}
$$

\begin{defi}

Soit $\epsilon$ un {\'e}l{\'e}ment de  $\mathcal{E}=\{0,1\}^N$ de longueur $l$
strictement positive. 
On note $\{i_1 < \cdots <i_l \}$ les {\'e}l{\'e}ments
 de $\pi_+(\epsilon)$. On d{\'e}finit alors l'application 
$r_M^{\epsilon} : \mathcal{U} \rightarrow \mathbb{Z}$ 
de la mani{\`e}re suivante :

\begin{enumerate}

\item[$(i)$] $r_M^{\epsilon}$ est $\mathbb{Z}$-lin{\'e}aire,


\item[$(ii)$] si $P \in \mathcal{U}$ est un mon{\^o}me non nul qui s'{\'e}crit 
sous la forme $P=SX_{i_l}^r Z_{i_l}^s$
  o{\`u} $S\in \mathcal{U}$ est de degr{\'e} $0$ en  $X_{i_l}$
   et $Z_{i_l}$, et o{\`u} $s$ est un entier positif et $r$
  un entier quelconque, 
$$r_M^{\epsilon}(P)=\left\{ \begin{array}{ll}
r_M^{\epsilon-(i_l)}\Big(S(1-m_{i_l})^{s-1}m_{i_l}^r\Big) & { \rm si \,
} s>0, \\
-r_M^{\epsilon-(i_l)}\Big(S(m_{i_l}+m_{i_l}^2 + \cdots + m_{i_l}^{r-1})\Big)
& { \rm si \,} s=0 { \rm  \, \, et \, } r>1, \\
r_M^{\epsilon-(i_l)}\Big(S(1+m_{i_l}^{-1} + \cdots
+ m_{i_l}^{r})\Big)
& { \rm si \,} s=0 { \rm \, \,  et \, } r < 0, \\
r_M^{\epsilon-(i_l)}(S) & { \rm si \,} s=0 { \rm \, \,
  et \, } r = 0, \\
0  & { \rm si \,} s=0 { \rm \, \,
  et \, } r = 1,
\end{array} \right.
$$
\item[$(iii)$] $r_M^{(\mathbf{0})}(P)=P(X_i=1, Z_i=0)$.

\end{enumerate}

\smallskip

Ces trois relations d{\'e}finissent compl{\`e}tement (r{\'e}cursivement) les applications
 $r_M^{\epsilon}$. On pose $r_M=r_M^{(\mathbf{1})}$.

\end{defi}

\smallskip

Pour $\epsilon \in \mathcal{E}$, on pose $s_{\epsilon}^T=\prod_{i \in
   \pi_-(\epsilon)} X_i \prod_{j \in \pi_+(\epsilon)}Z_j \in
   \mathcal{U}$. le th{\'e}or{\`e}me suivant est une cons{\'e}quence
   imm{\'e}diate du th{\'e}or{\`e}me~\ref{cskvd} :

\begin{theo} \label{cskovd}

Pour tout couple $(u,v) \in W^2$,
$$t_{u,v}^w = r_M\Big((\!\!\!\!\sum_{\epsilon \in \mathcal{E},
 \,
 \underline{v}(\epsilon)=\underline{u}}\!\!\!\!\!\!\!s_{\epsilon}^T)
 \,\,(\!\!\!\!\sum_{\epsilon' \in \mathcal{E},
 \,
 \underline{v}(\epsilon')=\underline{v}}\!\!\!\!\!\!\!s_{\epsilon'}^T) 
\Big).
$$

\end{theo}

\begin{exem}

On se place dans le cas $G_2$, et on calcule
$$t_{1,1}^{s_2s_1s_2s_1s_2}=r_M(X_1^2X_2^2X_3^2X_4^2X_5^2),$$
 o{\`u} $M=(\alpha_2,\alpha_1,\alpha_2,\alpha_1,\alpha_2)$. On calcule
 d'abord les mon{\^o}mes $m_i$ :
$$ \begin{array}{cll}
m_5 & = & X_1^{-2}X_2X_3^{-2}X_4,\\
m_4 & = & X_1^{3}X_2^{-2}X_3^3,\\
m_3 & = & X_1^{-2}X_2,\\
m_2 & = & X_1^3,               \\
m_1 & = & 1,
\end{array}
$$
ce qui nous donne 
$$t_{1,1}^{s_2s_1s_2s_1s_2}=-r_M^{(\overline{4})}(X_2^3X_4^3)
=r_M^{(\overline{3})}(X_1^3X_2X_3^3+X_1^6X_2^{-1}X_3^6) 
$$
$$=-r_M^{(\overline{2})}\big((X_1X_2^2+X_1^{-1}X_2^3)+(X_1^4+X_1^2X_2+
X_2^2+X_1^{-2}X_2^3+X_1^{-4}X_2^4)\big)
$$
$$=r_M^{(\overline{1})}\big(X_1^4+(X_1^2+X_1^5)-X_1^4+0+X_1^3+(X_1+X_1^4)+
(X_1^{-1}+X_1^2+ X_1^5)\big)
$$
$$=-3-1-4+3+0-2+0-3+2-1-4=-13.
$$

\end{exem}

\begin{rema}

On peut {\'e}galement restreindre les r{\'e}sultats des
sections~\ref{sectiontours} et \ref{sectionBott} au cas de la
$K$-th{\'e}orie ordinaire. On montre en effet facilement par
r{\'e}currence sur la dimension que la proposition~\ref{ordinaire} est
vraie pour les tours de Bott.

\end{rema}

   \bibliography{Ktheory}
   \bibliographystyle{smfplain}

  \end{document}